\documentclass{article}
\usepackage{amsmath} 
\usepackage{graphicx} 
\usepackage{mathtools} 

\DeclareMathOperator{\sgn}{sgn}

\title{Calculus:  a limitless perspective}
\author{Michael P. Lamoureux\thanks{University of Calgary (mikel@ucalgary.ca)}  
         $\,$ and Matt Yedlin\thanks{University of British Columbia (matty@ece.ubc.ca)} }
\date{October 14, 2025}

\begin{document}

\maketitle

\begin{abstract}
   We propose a novel foundation for calculus that focuses on the notion of approximations while avoiding the use of limits altogether. Continuity is defined as approximation at a point, while differentiability is defined as approximation with a linear function. The errors in  approximation are defined as a class of functions with certain properties; rules for combining error functions lead to all the familiar results in differential calculus. 

   We believe that this approach is more natural for students while still giving a rigourous foundation to differential calculus. We demonstrate its utility by deriving the basic differential rules for trigonometric, hyperbolic and exponential functions, as well as L'H\^{o}pital's Rule, Taylor polynomials, and the Fundamental Theorem of Calculus, all via approximation.
\end{abstract}

\section{Introduction}

The aim of this article is to demonstrate the development of a first-year calculus course without the use of limits. We focus on the notion of {\it approximations}, along with ideas from linear algebra and geometry, to give a solid foundation for a full course in calculus.

The reason for dropping limits is simple. Students find limits difficult in the introduction to freshman calculus, which unnecessarily increases their cognitive load and obscures the fundamental goals in calculus. In our experience, we find there are two principal reasons for the students' difficulties with limits. First,
their algebraic manipulation skills, so necessary in the discussion of limits,
are insufficient; the necessary elementary factorization patterns are not
embedded in their prerequisite backpack. Second, the necessary geometric background,
at a very basic level, is absent. Indeed, in many high school curricula, geometry is not
even presented. These deficiencies push the students towards memorizing formulas early in their
calculus journey, rather than understanding the basic principles of calculus. 

Yet limits seem to remain as the standard approach in our modern repertoire of calculus texts, such as in~\cite{apostol1967calculus, lang1986first, spivak2008calculus, stewart2015calculus, strang2019calculus, weir2017thomas}. In our opinion, limits are typically introduced as a procedural tool for students to
compute derivatives and integrals and possibly as a route to formal proofs.
These procedures also are convenient for instructors when creating tests for the students' work. 
Procedural methods and proofs can be important, but the use of limits in an
introductory course obscures the reason and power of doing calculus.

Developing calculus without the use of limits is not a new idea. For instance, Newton~\cite{newton1999principia} and Leibnitz both developed their pioneering methods using infinitesimals and differentials rather than limits. Building on the works of Amp\'{e}re~\cite{ampere1806recherches}, Euler~\cite{euler1755institutiones,euler1797introductio}, and Lagrange~\cite{lagrange1804leccons,lagrange1813theorie,lagrange1806lecons}, Cauchy created a rigourous foundation for calculus with the introduction of limits~\cite{cauchy1821cours, cauchy1826leccons,cauchy1829leccons}, but still mathematicians have looked for alternatives to the limit. For instance, the book by Marsden and Weinstein~\cite{marsden1981} uses the notion of a function ``overtaking'' its tangent as the basis for derivative, while the book by Zhang and Tong~\cite{zhang2018} introduce ``control functions'' that bound the difference ratio as the formulation for the derivative. Livshits~\cite{livshits2018simplify} pursues an algebraic approach while Sparks~\cite{sparks2004} has a more poetic viewpoint with his ``Calculus without limits - almost.'' Grabiner~\cite{grabiner1978,grabiner1983epsilon} presents an excellent review of the history of the development of the limit, as does Felscher~\cite{Felscher01112000} and Bair et al.~\cite{bair2013victors}.

Our novel contribution here is to focus on the notion of {\it approximation} as the formal way to define continuity and differentiability, and to {\it reify} the error term. That is, we {\it give a name} to the error function and define its properties that make the characteristics of the approximation precise. We specify what are the key properties of an error function, produce simple algebraic rules to combine error functions, and use these basic tools as a robust foundation for all the content of calculus. 

This approach drops limits altogether and appeals to the students' intuitive idea of an approximation. This seems like a simple idea for students to grasp. By adding some linear algebra and geometry, we have enough to form a rigorous basis for calculus. We demonstrate the power of this method by defining geometrically the trigonometric, hyperbolic and exponential functions and their derivatives, rather than the more traditional limit approach.\footnote{Such as $e^x = \lim_{n\to\infty} (1 + \frac{x}{n})^n$} We also use the method to derive an intuitive demonstration of  L'H\^{o}pital's rule as well as Taylor's theorem and the Fundamental Theorem of Calculus. 

\section{Introducing approximations}

We say two numbers are approximately equal if the difference between them is small. For instance, two numbers $a$ and $b$ are approximately equal if the difference $a-b$ is a small number. We often parameterize this statement by writing $a-b =\epsilon$, where $\epsilon = $ epsilon is some small  number, close to zero. Of course, this is a relative term. For example, we could say that two people are approximately the same height if the difference in their heights is less than one centimeter.  Two road trips are approximately the same distance if their distance is less than one kilometer. But not to worry, we will use this idea of an approximation as intuitive for now.

For functions, we say two functions are approximately equal if their difference is small, in some prescribed way. Our motivating examples are for continuity and for differentiability. For instance, we say a function $f(x)$ is {\em continuous} at a point $x_0$ if it is approximately equal to the constant $f(x_0)$ for all $x$ near $x_0$. We say that a function $f(x)$ is {\em differentiable} at a point $x_0$ if it is approximately equal to the constant $f(x_0)$ plus a linear term $L(x-x_0)$ for all $x$ near $x_0$.   These motivating examples extend naturally to higher dimensions and more complicated spaces without any mention of limits or partial derivatives, facilitating the students' understanding of both introductory and advanced calculus.

We do need to be careful by what we mean by ``{\em approximate}'' for functions. In the case of continuity, the difference between the function value $f(x_0 + \epsilon)$ and the constant $f(x_0)$ is some small error function. In the derivative case, we will mean the  difference between the function value $f(x_0 + \epsilon)$ and the linearized approximation $f(x_0) + L(\epsilon)$ is equal to some other error function times the size of $\epsilon$. We write out in full as 
\begin{eqnarray*}
    f(x_0 + \epsilon) &=&  f(x_0) + E_0(\epsilon), \mbox{ for continuity}, \\
    f(x_0 + \epsilon) &=&  f(x_0) +  L(\epsilon) +  |\epsilon|E_1(\epsilon), \mbox{ for differentiability} .
\end{eqnarray*}
These error functions $E_0(\epsilon),E_1(\epsilon)$ can depend on both $f$ and $x_0$ but must have these three key properties:
\begin{enumerate}
    \item $E_*(\epsilon)$ is defined for all small values of $\epsilon$;
    \item $E_*(0) = 0$;
    \item $E_*(\epsilon)$ is made small by making $\epsilon$ small.
\end{enumerate}
This third statement needs to be made precise without appealing to limits, which we will address later. We are of course thinking of error functions in the form like $E(\epsilon) = \epsilon^2$ or $E(\epsilon) = |\epsilon|^{1/2}$ or similar functions that are bounded by monomials in $\epsilon$ or its roots.

In this motivating example, we say the linear map $L(\epsilon)$ is the derivative of $f(x)$ at the point $x=x_0$ and we write $L(\epsilon) = f'(x_0)\cdot\epsilon$, where $f'(x_0)$ is defined the scaling constant for the linear map $L(\epsilon)$. The approximation is then expressed as
$$f(x_0 + \epsilon) =  f(x_0) +  f'(x_0)\cdot\epsilon +  |\epsilon|E(\epsilon).$$  With this example as motivation\footnote{We use the absolute value $|\epsilon|$ in the definition, looking ahead to multivariable calculus. This also agrees with the definition of the Fr\'{e}chet derivative~\cite{cartan1971}.} we will build up calculus without using limits. 

Note the use of error functions extends naturally to other concepts in calculus. For example, a sequence of numbers $a_n$ converges to the point $a$ if there is an approximation of the form \[ a = a_n + E(\frac{1}{n})\] for some error function $E(\frac{1}{n})$ which goes to zero as $n$ increases without bound. Similarly, a sequence of functions $f_n(x)$ converges pointwise to $f(x)$ if there are error functions $E_x(\epsilon)$ at each point $x$ yielding the approximations \[ f(x) = f_n(x) + E_x(\frac{1}{n}). \]  A Taylor expansion for differentiable function $f(x)$ is a polynomial $p(x)$ of order $n$ with the approximation identity 
\[ f(x_0 + \epsilon) = p(\epsilon) + |\epsilon|^n E(\epsilon)\]
for some error function $E(\epsilon)$. And so on. Of course, this is very much along the lines of the ``little-o'' order of approximation for functions~\cite{apostol1967calculus}. However, the important innovation here is we give a name to the error as a function, in order to clarify its properties and manipulate it directly in calculations.

\section{The algebra of error functions}

We make the formal definition of an error function precise with the following statement. 

{\bf Definition:} A function $E(\epsilon)$ is called an {\bf error function} if 
\begin{enumerate}
    \item $E(\epsilon)$ is defined for all small values of $\epsilon$;
    \item $E(0) = 0$;
    \item $E(\epsilon)$ can always be made as small as desired by making $\epsilon$ sufficiently small.
\end{enumerate}

Property 3) is key to the idea of a good approximation and merits further explanation. It is not just that the error $E(\epsilon)$ is small when $\epsilon$ is small, but in fact, if some tolerance is specified, such as $-0.15 < E(\epsilon) < 0.15$, then one can always find some small numbers $x_0<0<x_1$ so that $\epsilon$ between $x_0$ and $x_1$ will imply $E(\epsilon)$ is within that tolerance $\pm 0.15$. The important point is that finding such numbers $x_0,x_1$ is possible for any tolerance specified, not just $\pm 0.15$. 

\begin{figure}
    \centering
    \includegraphics[width=0.7\linewidth]{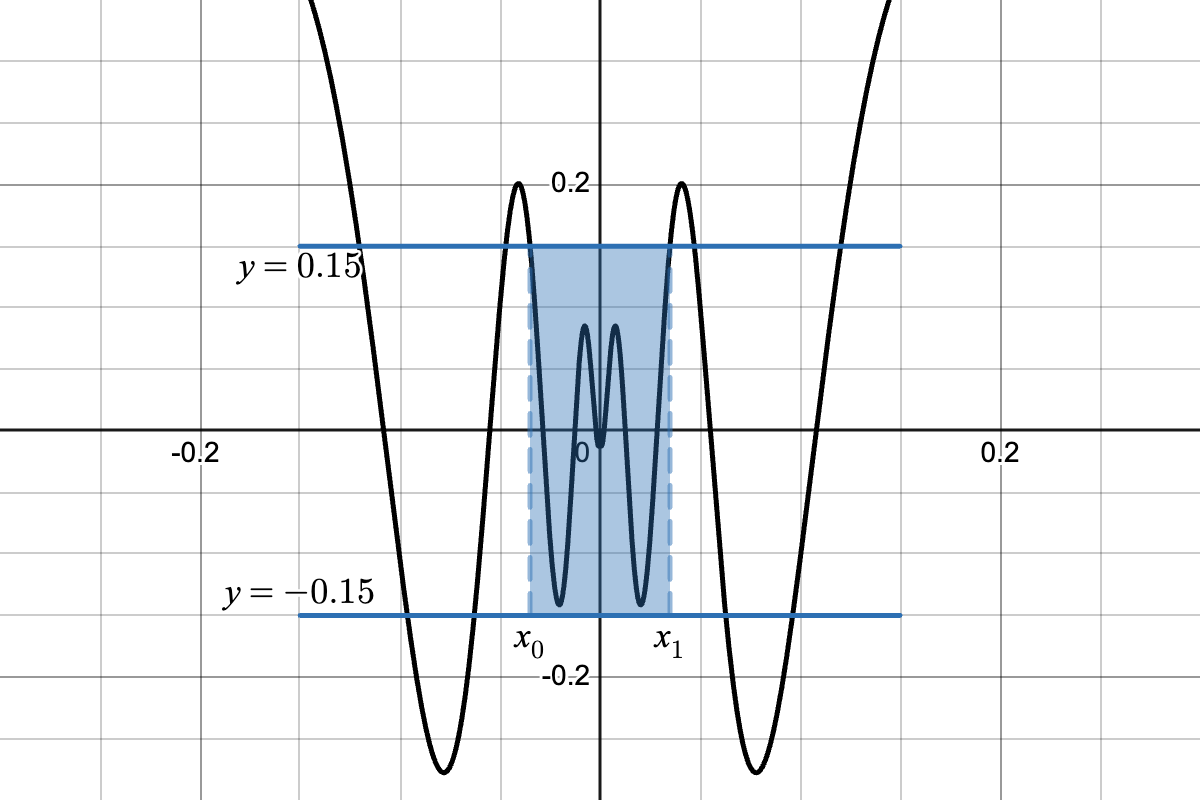}
    \caption{Bounding box on an error function, demonstrating Property 3.}
    \label{fig:ErrorFunction}
\end{figure}

For the beginner, this may be enough to get the intuitive idea of an error function. A more formal definition is demonstrated by Figure~\ref{fig:ErrorFunction}, and stated as follows:

{\bf Definition: } Property 3 means that, given any specified tolerance $y_0<0<y_1$ there is an interval $[x_0,x_1]$, with $x_0<0<x_1$ so that $\epsilon$ in $[x_0,x_1]$ implies $E(\epsilon)$ is in the interval $[y_0,y_1]$.

In other words, the graph of the function $E(\epsilon)$, restricted to the interval $[x_0,x_1]$, sits inside the rectangle $[x_0,x_1]\times [y_0,y_1]$, as shown in Figure~\ref{fig:ErrorFunction}. 

\begin{figure}
    \centering
    \includegraphics[width=0.7\linewidth]{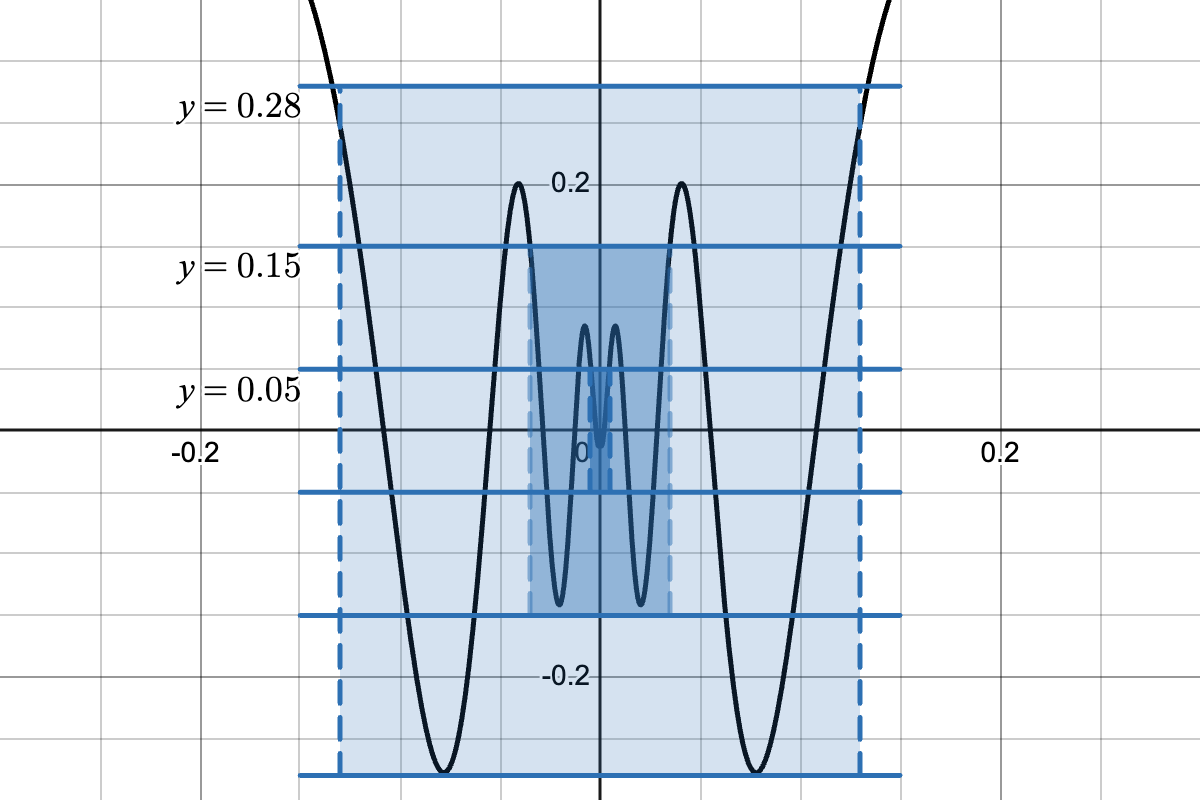}
    \caption{A sequence of bounding boxes, demonstrating Property 3.}
    \label{fig:ErrorFunction_3}
\end{figure}

\begin{figure}
    \centering
    \includegraphics[width=0.7\linewidth]{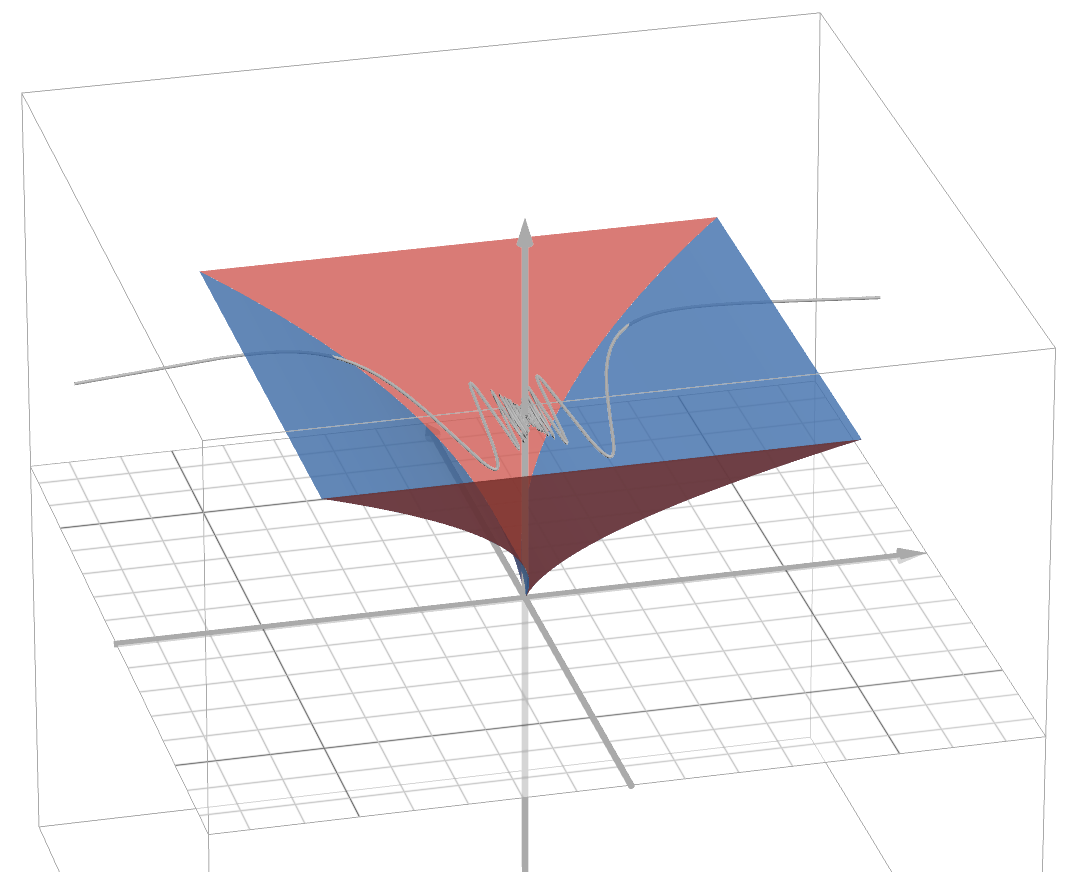}
    \caption{A funnel to visualize property 3 of the error function.}
    \label{fig:ErrorFunction_funnel}
\end{figure}

In fact, we should think of a nested sequence of rectangles telescoping inwards towards the origin, as shown in in Figure~\ref{fig:ErrorFunction_3}. As the height of the rectangles decreases to zero, the width of the rectangles decreases to zero fast enough to always capture the graph of the function within the rectangle. It is instructive to visualize this sequence of rectangles as a 3D funnel shrinking to the origin, as shown in Figure~\ref{fig:ErrorFunction_funnel}. Again, as the funnel shrinks in the $y$-direction, the $x$-direction must also shrink fast enough to contain the graph of the error function. 

Some examples of error functions are $E(\epsilon) = k\epsilon$ for some constant $k$, $E(\epsilon) = |\epsilon|^{1/2}$, or $E(\epsilon) = \frac{\sin(\epsilon)-\epsilon}{|\epsilon|}$.

As a consequence of the definition, we have the following algebraic properties in combining error functions:
\begin{itemize}
    \item The sum and difference of two or more error functions is an error function;
    \item The product of an error function times a bounded function is again an error function;
    \item The composition of one error function with another error function is again an error function;
    \item if an error function $E(\epsilon)$ bounds another function $E_1(\epsilon)$, which is to say $|E(\epsilon)| \geq |E_1(\epsilon)|$ for all small $\epsilon$, then $E_1(\epsilon)$ is also an error function;
    \item An error function $E(\epsilon)$ is bounded for $\epsilon$ sufficiently small;
    \item For any number $M>0$, an error function $E(\epsilon)$ is bounded above and below by $\pm M$ for $\epsilon$ sufficiently small.
\end{itemize}

These properties will seem natural to students: for instance, the sum rule just says if we have two small errors, they add up to a small error. Demonstrating these properties is straightforward, although the level of argument may depend on the mathematical level of students learning in class. 

For example, to show the sum property, if $E_1(\epsilon),E_2(\epsilon)$ are two error functions, then the sum $E(\epsilon) = E_1(\epsilon) + E_2(\epsilon)$ is also an error function since:
\begin{itemize} 
\item Property 1) is satisfied, as the summands are defined for all small $\epsilon$ and thus the sum is defined;
\item Property 2) is satisfied as $E(0) = E_1(0) + E_2(0) = 0+0 = 0$;
\item Property 3) is satisfied as, if we take interval $[x_0,x_1]$ small enough that both $E_1([x_0,x_1])$ and $E_2([x_0,x_1])$ lie in the interval $[y_0/2,y_1/2]$, then the sum $E_1 + E_2$ lies in the interval $[y_0/2,y_1/2]+[y_0/2,y_1/2] = [y_0,y_1]$.
\end{itemize}

\section{Definition of continuity}

To define continuity, we make the statement that a function $f(x)$ is continuous at a point $x=x_0$ if it is well-approximated by the constant function $f(x_0)$. Or in other words, if $x$ is close to $x_0$, then $f(x)$ is close to $f(x_0)$. We make this statement precise by stating the difference $f(x_0 + \epsilon) - f(x_0)$ is given by some error function $E(\epsilon)$ for all small $\epsilon$. Writing this out in full, we say $f$ is continuous at $x_0$ if $$f(x_0 + \epsilon) =  f(x_0) + E(\epsilon), \mbox{ for $\epsilon$ small,} $$ where $E(\epsilon)$ is one of these error functions as defined in Section~2. 

The attentive reader will note that any function $f(x)$ that is defined for values of $x$ near 0, and is continuous at $x=0$ with value $f(0) = 0$, then it is in fact an instance of an error function as defined in Section~2. And conversely. 

A function $f$ is continuous on a set if it is continuous at every point in that set.

Some examples of continuous functions are a constant function (like $f(x) = 10$), a linear function (like $f(x) = 10x$) and a polynomial function (like $f(x) = 10x^3 - 9x^2 + 8x -7$). In these cases, suitable  error functions that show continuity can be determined explicitly. 

\subsection{Combining continuous functions}

Given two continuous functions $f,g$, the following combinations are also continuous:
\begin{itemize}
    \item The sum and differences $f(x)+g(x)$, $f(x) - g(x)$ are continuous;
    \item The product $f(x)g(x)$ is continuous;
    \item The reciprocal $1/g(x)$ is continuous at    all points where $g(x)$ is not zero;
    \item The quotient $f(x)/g(x)$ is continuous at all points where $g(x)$ is not zero;
    \item The composition $f(g(x))$ is continuous.
\end{itemize} 

These results follow immediately from the properties of the error functions. 

For example, to demonstrate the sum rule, fix a point $x_0$ where $f, g$ are continuous, so we have the approximations
$$f(x_0 + \epsilon) = f(x_0) + E_1(\epsilon), \qquad
g(x_0 + \epsilon) = g(x_0) + E_2(\epsilon).$$
Then the sum $h(x) = f(x)+g(x)$ satisfies
$$h(x_0 + \epsilon) = f(x_0) + E_1(\epsilon) + g(x_0) + E_2(\epsilon) = h(x_0) + (E_1(\epsilon) + E_2(\epsilon)),$$
so $h(x_0+\epsilon)$ is given as $h(x_0)$ plus the error function $E(\epsilon) = E_1(\epsilon) + E_2(\epsilon).$ Thus $h$ is continuous at $x=x_0$ by definition.

To show the reciprocal rule, we first show the function $f(x) = 1/x$ is continuous at all points $x_0\neq 0.$ We write
\begin{eqnarray*}
    f(x_0 + \epsilon) &&=
    \frac{1}{x_0 + \epsilon} \\ &&=
    \frac{1}{x_0} - \frac{1}{x_0} + \frac{1}{x_0 + \epsilon}  \\
    &&=
    \frac{1}{x_0} + \frac{x_0 - (x_0 + \epsilon)}{x_0(x_0+\epsilon)} \\
    &&= 
    \frac{1}{x_0} - \epsilon \frac{1}{x_0(x_0+\epsilon)} \\
    &&= 
    f(x_0) + E(\epsilon), \\
\end{eqnarray*}
where $E(\epsilon) = \epsilon/(x_0(x_0 + \epsilon))$ is an error function, since it is just $\epsilon$ times a bounded function (for $\epsilon$ small). This shows $f$ is continuous. Now consider the composition $f(g(x)) = 1/g(x)$. By the composition rule for continuous functions, since $f$ and $g$ are both continuous, this composition is also continuous, at the points where the denominator $g(x)$ is non-zero. 

The demonstration of the other combination rules is left as an exercise.

\section{Definition of derivative}

We say a function $f(x)$ is differentiable at a point $x=x_0$ if there exists a linear function $L(x)$ which yields a first-order approximation to $f(x)$ for points near $x_0$ which is also the best first-order approximation. More precisely we say the function $f(x)$ is differentiable at $x_0$ if this equality  $$f(x_0 + \epsilon) = f(x_0) + L(\epsilon) + |\epsilon|E(\epsilon)$$
holds for all $\epsilon$ sufficiently small, where $E(\epsilon)$ is an error function, as defined in Section~3. 

This choice of error in the form of $|\epsilon|E(\epsilon)$ ensures that the linear approximation determined by $L$ is a good approximation, and is also the unique, optimal choice, as we will see in the next section.

When such an approximation holds, we say the linear map $L()$ is the derivative of $f(x)$ at the point $x=x_0$. In this case, we write the linear map $L()$ as multiplication by a constant $a$, which defines the derivative as $f'(x_0)=a$ and the approximation is given as 
$$f(x_0 + \epsilon) = f(x_0) + f'(x_0)\cdot\epsilon + |\epsilon|E(\epsilon).$$

The choice of absolute value $|\epsilon|$ in the approximation is deliberate, as it simplifies certain algebraic manipulations in the work that follows. It also generalizes easily to multivariate calculus. 

\section{The derivative is unique}

The derivative of a function will be unique. To see this directly, suppose we have two candidates for a first order approximation for the function. We can write
\begin{eqnarray*}
    f(x_0 + \epsilon) &= f(x_0) + a_1\cdot \epsilon + |\epsilon| E_1(\epsilon), \\
    &= f(x_0) + a_2\cdot \epsilon + |\epsilon| E_2(\epsilon).
\end{eqnarray*}
Taking the difference gives the equation
$$0 = (a_1 - a_2)\cdot\epsilon + |\epsilon|(E_1(\epsilon) - E_2(\epsilon))$$
and dividing by $\epsilon $ gives
$$a_2 - a_1 = \sgn(\epsilon)(E_1(\epsilon) - E_2(\epsilon)), $$
where $\sgn(\epsilon)$ is the sign function. Letting $\epsilon$ tend to zero puts a zero on the right hand side, so we have $a_2-a_1 = 0$. Hence the two possible values for the derivative are in fact the same.

\section{Derivative of a monomial}

Let us compute some simple derivatives using the approximation definition. 

As a first example, consider the function $f(x) = x^2$. Set $x_0 = 3$ so $f(x_0) = 3^2 = 9$. Then $f(x_0 + \epsilon) = (3 + \epsilon)^2 = 9 + 6 \epsilon + \epsilon^2$. We have the identity $$f(x_0+ \epsilon) = f(x_0) + 6\epsilon + \epsilon^2 = f(x_0) + 6\epsilon + |\epsilon| E(\epsilon),$$ where we have the error function given as $E(\epsilon) = |\epsilon|$. This shows the function $f(x) = x^2$ has a first order linear approximation at $x=3$ and thus it is differentiable there, with derivative $f'(3) = 6$. 

In general, for a monomial function of the form $f(x) = x^n$, with $n>2$, we have $$f(x+\epsilon) = (x+\epsilon)^n = \sum_{k=0}^n \binom{n}{k} x^{n-k} \epsilon^k = x^n + (nx^{n-1})\epsilon + \sum_{k=2}^n \binom{n}{k}x^{n-k} \epsilon^k,$$ where we used  the binomial expansion and pulled out the first two terms. The last term, with the summation, has factors of $\epsilon^2$ or higher so we can express it as in the form of $|\epsilon|$ times an error function to obtain
$$f(x+\epsilon) = f(x) +  (n x^{n-1}) \epsilon + |\epsilon|\left[ |\epsilon| \sum_{k=2}^n \binom{n}{k}x^{n-k} \epsilon^{k-2}\right].$$
Here, the error function is just the term in the square brackets, \linebreak $E(\epsilon) = |\epsilon| \sum_{k=2}^n \binom{n}{k}x^{n-k} \epsilon^{k-2}.$
(Why is this an error function? Well, for $x$ fixed, it is a sum and product of continuous functions in $\epsilon$ and thus is continuous, and take the value zero at $\epsilon = 0$. So yes, that is an error function.)

Thus from the approximation formula, the linear term in $\epsilon$ shows the monomial $f(x) = x^n$ is differentiable at any $x$, with derivative $f'(x) = nx^{n-1}$. 

\section{Derivative rules}

The approximation formulas and properties of the set of error functions makes it easy to show the familiar differentiation rules for the sum, difference, product and quotients of  differentiable function. The chain rule is also a simple consequence. 

We list the rules here. Suppose $f,g$ are differentiable functions. Then we have the following rules for combining derivatives:
\begin{itemize}
    \item Sum rule: $$(f + g)'(x) = f'(x) + g'(x)$$
    \item Difference rule: $$(f - g)'(x) = f'(x) - g'(x)$$
    \item Product rule: $$(f\cdot g)'(x) = f'(x)g(x) + f(x)g'(x)$$
    \item Quotient rule: $$(\frac{f}{g})'(x) = \frac{f'(x)g(x) - f(x)g'(x)}{(g(x))^2}$$
    \item Chain rule: $$(f\circ g)'(x) = f'(g(x))\cdot g'(x))$$
    \item Inverse function rule: For $g = f^{-1}$, an inverse function, then
    $$g'(x) = \frac{1}{f'(g(x))}$$
\end{itemize}

As an example, let us verify the product rule, to show the utility of the approximation approach. It nicely shows why a sum of two products arise in the derivative of the product. 

With $f,g$ differentiable, we have the approximations
\begin{eqnarray*}
    f(x + \epsilon) &&= f(x) + f'(x)\cdot\epsilon + |\epsilon|E_1(\epsilon), \\
    g(x + \epsilon) &&= g(x) + g'(x)\cdot\epsilon + |\epsilon|E_2(\epsilon). 
\end{eqnarray*}
Thus, the product of the approximations satisfies
\begin{eqnarray*}
    (f\cdot g)(x + \epsilon) &&= 
    [f(x+\epsilon)]\cdot [g(x + \epsilon] \qquad \qquad \qquad \qquad \qquad \qquad.\\
    &&= 
    [f(x) + f'(x)\cdot\epsilon + |\epsilon|E_1(\epsilon)]\cdot 
    [g(x) + g'(x)\cdot\epsilon + |\epsilon|E_2(\epsilon)] \\
    &&= 
    f(x)g(x) + [f'(x)g(x) + f(x)g'(x)]\cdot\epsilon + |\epsilon|E(\epsilon)]
\end{eqnarray*}
where $E(\epsilon) = E_1(\epsilon)(g(x + \epsilon) + f(x + \epsilon)E_2(\epsilon)$ is the required error function. The middle term in the formula, the part that is linear in $\epsilon$, shows the derivative of the product is $f'(x)g(x) + f(x)g'(x)$.

A demonstration of the chain rule is similarly straightforward. With $x$ fixed and $y = g(x)$,  we have the approximations
\begin{eqnarray*}
    f(y + \epsilon) &&= f(y) + f'(y)\cdot\epsilon + |\epsilon|E_1(\epsilon), \\
    g(x + \epsilon) &&= g(x) + g'(x)\cdot\epsilon + |\epsilon|E_2(\epsilon). 
\end{eqnarray*}
Setting $\tilde{\epsilon} = g'(x)\epsilon + |\epsilon|E_2(\epsilon)$ we $g(x+\epsilon) = g(x) + \tilde{\epsilon} = y +\tilde{\epsilon}$. Using the linear approximation with $f$ at the point $y = g(x)$ we obtain
\begin{eqnarray*}
    f(g(x+\epsilon)) && = f(y + \tilde{\epsilon}) \\
    && = f(y) + f'(y)\cdot\tilde{\epsilon} + |\tilde{\epsilon}|E_1(\tilde{\epsilon}) \\
    && = f(y) + f'(y))\cdot[g'(x)\cdot\epsilon + |\epsilon|E_2(\epsilon)] + |\tilde{\epsilon}|E_1(\tilde{\epsilon}) \\
    && = f(g(x)) + [f'(y)g'(x)]\cdot\epsilon + |\epsilon|E(\epsilon),
\end{eqnarray*}
where we have gathered all the error terms into the last term of this equation. From this equation, we see the middle term is the linear part with respect to $\epsilon$ with coefficient $f'(y)g'(x)$. and so we have the derivative of the composition $f\circ g$ is $(f\circ g)'(x) =f'(y)g'(x),$ where $y = g(x).$

The inverse function rule for derivatives follows immediately from the chain rule. For if $g$ is the inverse function for $f$, then $(f\circ g)(x) = x$ and differentiate both sides of the equation to obtain, via the chain rule,
\[ f'(y)g'(x) = 1, \mbox{ where } y=g(x).\]
Dividing out the derivative gives
\[ g'(x) = \frac{1}{f'(y)}, \mbox{ or more simply } g'(x) = \frac{1}{f'(g(x))}.\]

\section{Chain rule, visually}

\begin{figure}
    \centering
    \includegraphics[width=0.7\linewidth]{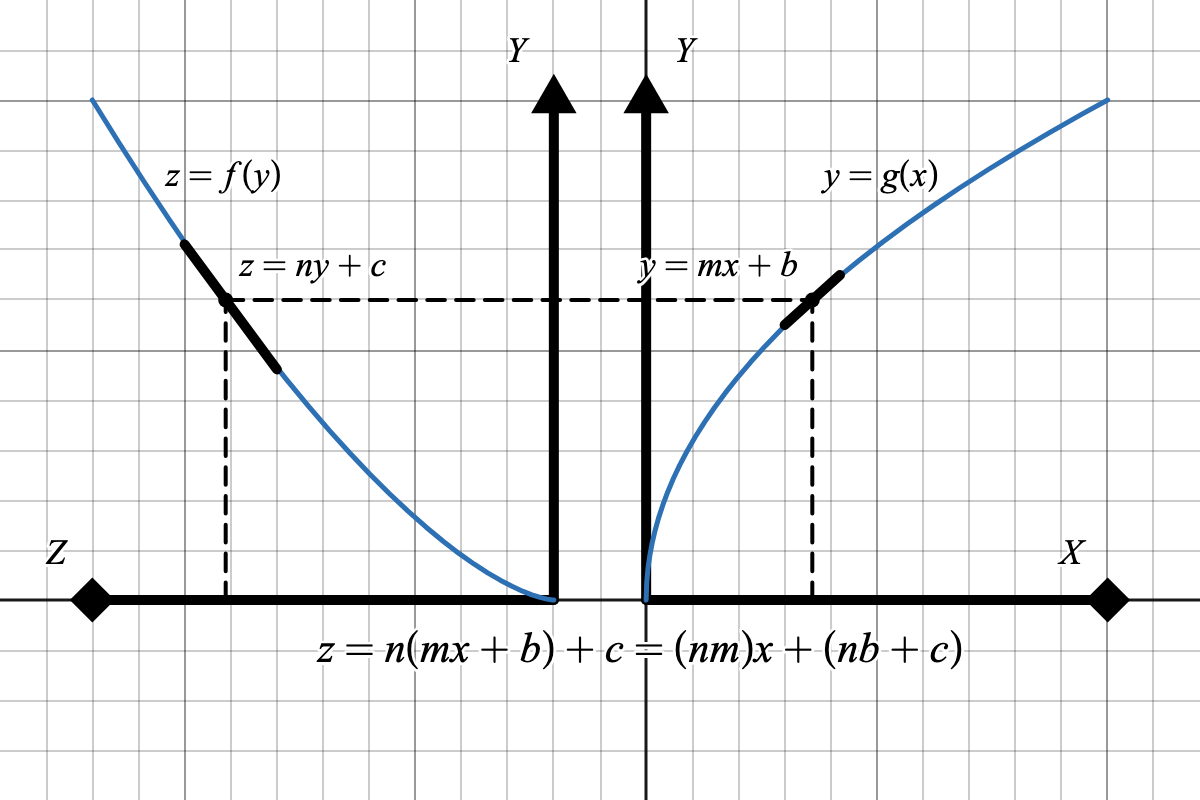}
    \caption{Composing two functions and the chain rule.}
    \label{fig:ChainRule}
\end{figure}

The chain rule can be easily visualized by putting together plots of the two functions $y=g(x)$ and $z=g(y)$, while rotating the second graph by 90 degrees to line up the two $y$-axes, as shown in Figure~\ref{fig:ChainRule}. The plot on the right shows $y=g(x)$ and the plot on the left shows $z=f(y).$ The linear approximations for the functions are given as $y = mx + b$ and $z = ny +c$. Composing the two functions leads to a composition of the approximations, so
\[ z = ny+c = n(mx+b) + c = (mn)x + (nb+c).\]
Thus, with $m$ the derivative of $g$ at a point $x_0$, $n$ the derivative of $f$ at the corresponding point $y_0$, the composition is a linear function with linear coefficient $mn$, which is the product of the two derivatives. 

\section{Rational powers of x}

Finding the derivative of a rational power of $x$ is straightforward using the chain rule. Suppose $f(x) = x^r$ for a rational number $r = m/n$ with positive integers $m,n$ and variable $x>0$. Then \[(f(x))^n = x^m,\]
and differentiating both sides, using the chain rule and the results on monomials,  gives
\[ n(f(x)))^{n-1}f'(x) = mx^{m-1}.\]
Thus
\[ f'(x) = \frac{m}{n}\frac{x^{m-1}}{(x^{m/n})^{n-1}}= r\frac{x^{m-1}}{x^m x^{-r}} = r x^{r - 1}.\]
Thus, for any positive rational $r$, with $f(x) =x^r$, we have the derivative \[ f'(x) = rx^{r-1}.\] 

The same result holds for negative $r$, which we can verify using the quotient rule for derivatives. 

\section{Demonstration: Calculus from geometry}

In the following three sections, we demonstrate how our approach of approximation and error functions can be applied directly to the construction of a geometric definition for the standard trigonometric, hyperbolic, and exponential functions. Without limits, we determine their derivatives as well as their basic algebraic properties. 

For these three classes of functions, the key geometric idea is the same. Namely, a quadratic equation defines a curve in the $xy$-plane and two function values $f(A),g(A)$ are determined as the first and second components of the unique point on that curve that selects a region of the given area of size $A$ or $A/2$. 

\begin{figure}
    \centering
    \includegraphics[width=0.45\linewidth]{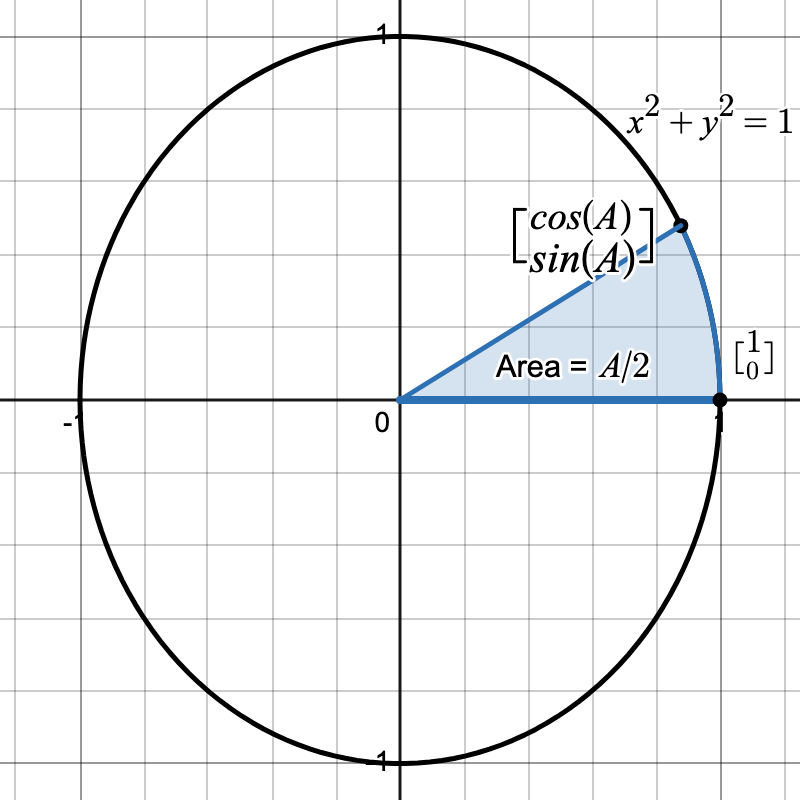}
    \includegraphics[width=0.45\linewidth]{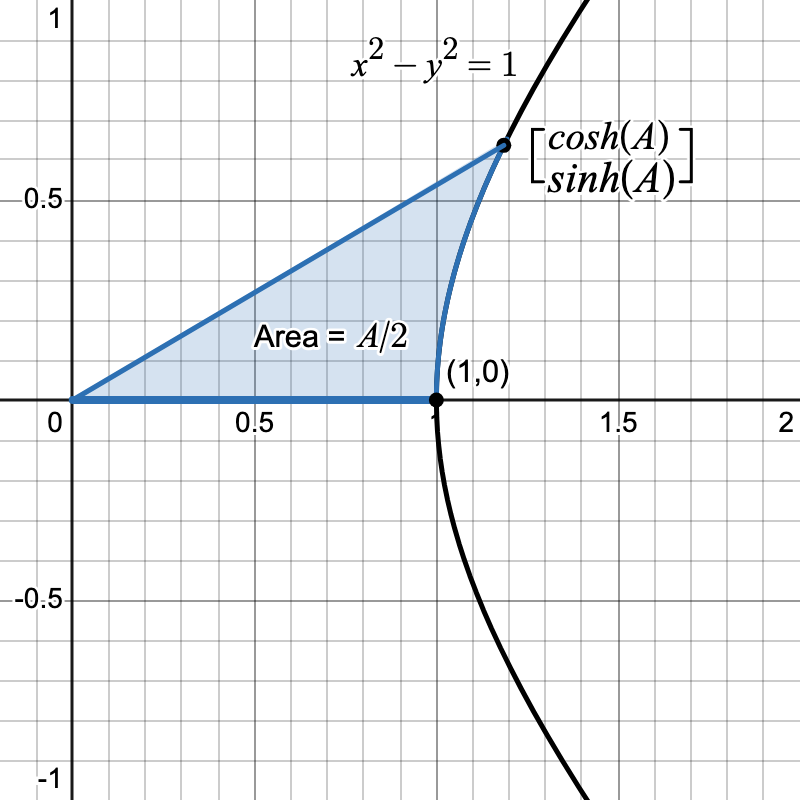}
    \includegraphics[width=0.45\linewidth]{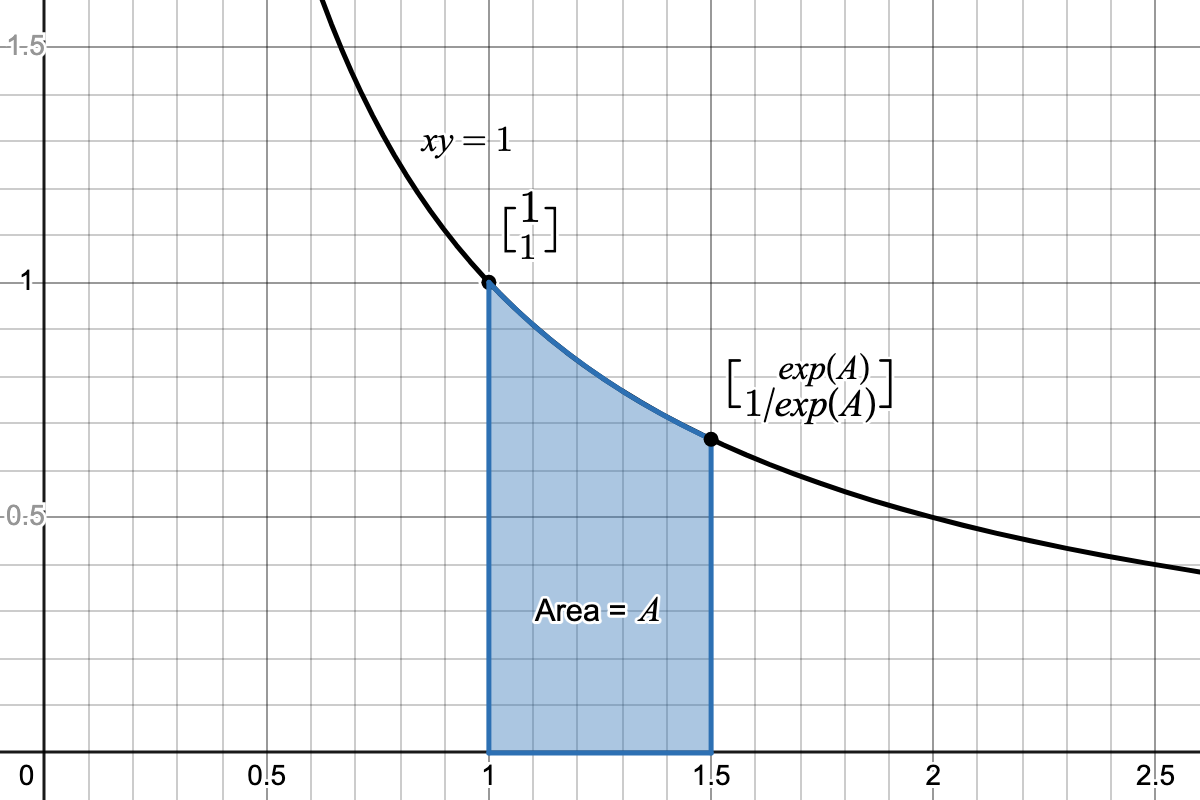}
    \caption{Three curves $x^2 + y^2 = 1$, $x^2-y^2 = 1$, $xy=1$, and associated regions defining trigonometric, hyperbolic, and exponential functions.}
    \label{fig:Function_defn}
\end{figure}

We use column notation 
$\begin{bsmallmatrix} x \\y \end{bsmallmatrix}$
for points in the plane, in order to facilitate the use of linear algebra in the forthcoming discussions. 

In Figure~\ref{fig:Function_defn} we see the three curves (a unit circle, a unit hyperbola and a skew rectangular hyperbola) and the related regions of given area. From the geometry, it is apparent that for small values of area $A$, one can always find a point 
$\begin{bsmallmatrix} x \\y \end{bsmallmatrix}$
on the given curve that defines a region of the required area.\footnote{In fact knowing that this point exists requires both continuity, and completeness property of the real numbers. But this is usually beyond what we teach students in first year.} Also, small changes in $A$ give only a small change in the selected region and thus only a small change in the point 
$\begin{bsmallmatrix} x \\y \end{bsmallmatrix}$,
and conversely, showing continuity of the functions $f(A),g(A)$ so defined. 

In the following sections, we will focus on the careful definition of these functions, their extensions to large values for the area $A$, their algebraic properties, and their derivatives. In each case, an area-preserving transformation is introduced to manipulate the areas under the curves.

\section{The trigonometric functions}

In this section, we define the familiar cosine and sine functions that arise in the geometry of circles and triangles in the plane. Rather than angles and arclength, we use {\it areas} of sections of a circle. This approach generalizes to the hyperbolic and exponential functions, shown in the next two sections, and avoids the problematic issue of determining the arclength of a curve. 

We will consider the basic definitions for smaller parameter values, then extend them to arbitrarily large values, demonstrate the summation rules, and finally compute the derivatives.

\subsection{Defining cosine and sine for small arguments}

\begin{figure}
    \centering
    \includegraphics[width=0.45\linewidth]{4_Circ_defn_pos.png}
    \includegraphics[width=0.45\linewidth]{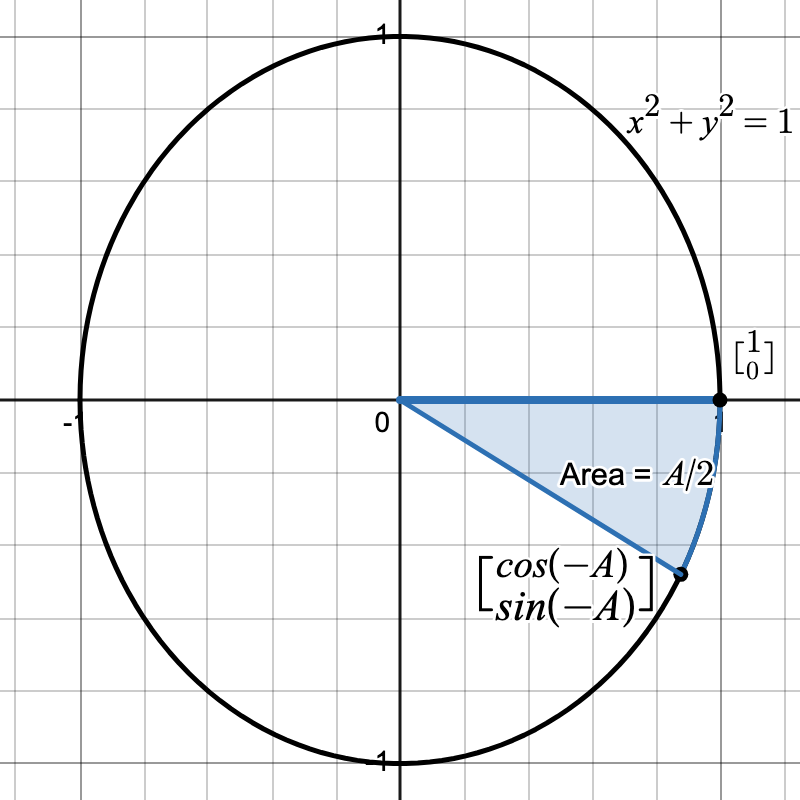}
    \caption{Given area $A/2$, functions $\cos(A), \sin(A)$ are defined.}
    \label{fig:Circ_defn}
\end{figure}

We begin with the unit circle $x^2 + y^2 =1$ on the plane, which we know from geometry has a total area of  $\pi\approx 3.14$. Given a positive number $A$, say smaller than $\pi$, we can find a circular section of the circle sitting above the positive $x$-axis whose area is exactly $A/2$. This section has a top vertex at a position $\begin{bsmallmatrix}x \\ y \end{bsmallmatrix}$ on the unit circle, which we use to define the values $\begin{bsmallmatrix}\cos(A)=x\\ \sin(A)=y \end{bsmallmatrix}$ as shown on the left in Figure~\ref{fig:Circ_defn}.

For the corresponding negative value $-A$, we drop a circular section below the $x$-axis, with area $A/2$, as shown in the right side of Figure~\ref{fig:Circ_defn}. This section has a lower vertex at a point $(x,y)$ below the horizontal axis, which we use to find the cosine and sine as 
$\begin{bsmallmatrix}\cos(-A)=x\\ \sin(-A)=y \end{bsmallmatrix}$.

 From the symmetry in Figure~\ref{fig:Circ_defn} we observe that including a negative sign in the argument does not change the value of the cosine, while it introduces a negative sign on the sine. Since the points are on the unit circle, we observe that the sum of the squares of the cosine and sine is one. We have obtained three key trigonometric relationships:
\begin{align*}
    \cos(-A) &= \cos(A) \\
    \sin(-A) &= -\sin(A) \\
    \cos^2(A) + \sin^2(A) &= 1 
\end{align*}

From simple geometry, we can evaluate the trig functions at some typical values:

\begin{align*}
    \cos(0) =&1,  &\sin(0) =& 0  \\
    \cos(\pi/4) =&1/\sqrt{2}, &\sin(\pi/4) =& 1/\sqrt{2}  \\
    \cos(\pi/2)=&0, &\sin(\pi/2) =& 1  \\
    \cos(3\pi/4) =&-1\sqrt{2},  & \quad\sin(3\pi/4) =& 1/\sqrt{2}  \\
    \cos(\pi) =&-1,  &\sin(\pi) =& 0   
\end{align*}

\subsection{Extending to large arguments}

\begin{figure}
    \centering
    \includegraphics[width=0.7\linewidth]{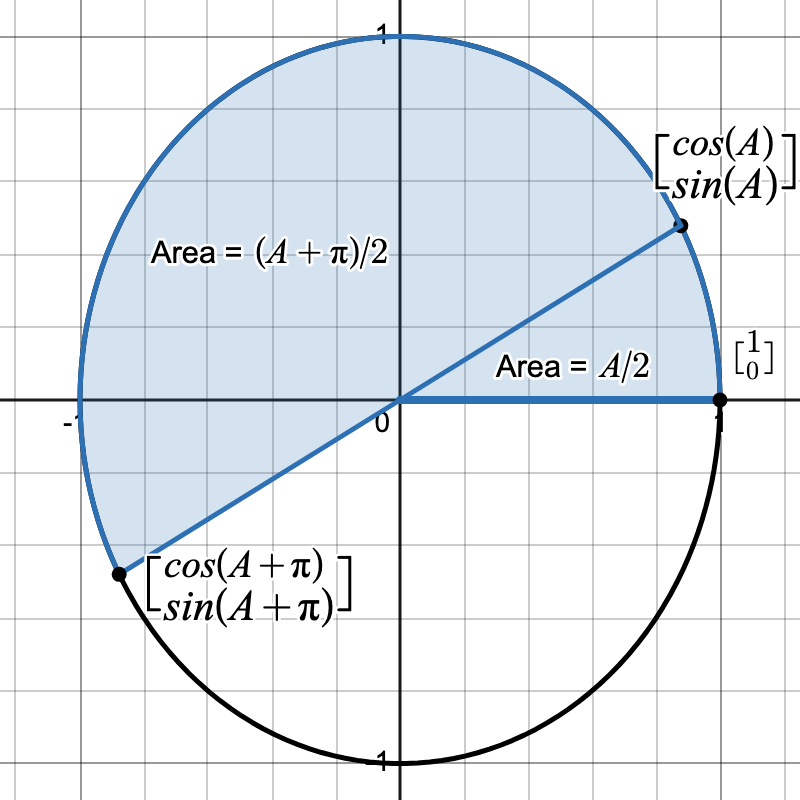}
    \caption{Adding an area of $\pi/2$ to define $\cos(A+\pi),\sin(A+\pi)$.}
    \label{fig:Circ_plus_pi}
\end{figure}

For larger arguments $A$, we will extend the functions $\cos()$, $\sin()$ by periodicity. Referring to Figure~\ref{fig:Circ_plus_pi}, we evaluate the functions at argument $A + \pi$ by adding a half circle of area $\pi/2$ to find that 
\begin{align*}
    \cos(A + \pi) =&-\cos(A),  &\sin(A+\pi) =& -\sin(A).  \\
\end{align*}
By adding multiples of $\pi$, we can extend the definition of the trig functions to arbitrarily large arguments, using the formulas: 
\begin{align*}
    \cos(A + n\pi) =&(-1)^n \cos(A),  &\sin(A+n\pi) =& (-1)^n\sin(A),  \\
\end{align*} 
for any integer $n$. Thus, we have the trig functions defined on the whole real line. 

\subsection{Summation rules}

In this section, we derive the familiar summation rules for determining the values of $\cos(A+B)$ and $\sin(A+B)$. 

\begin{figure}
    \centering
    \includegraphics[width=0.49\linewidth]{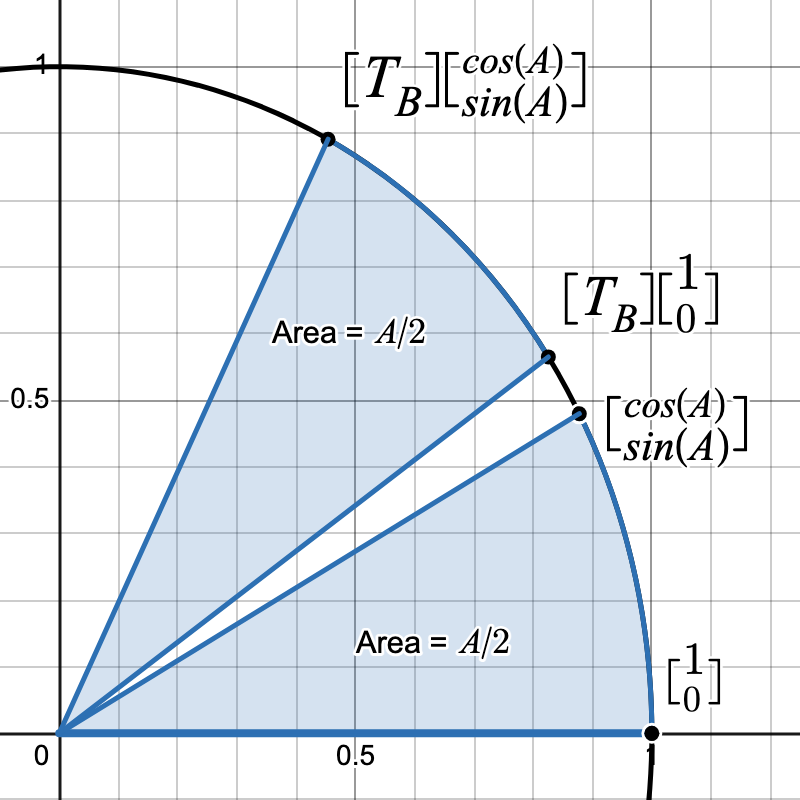}
    \includegraphics[width=0.49\linewidth]{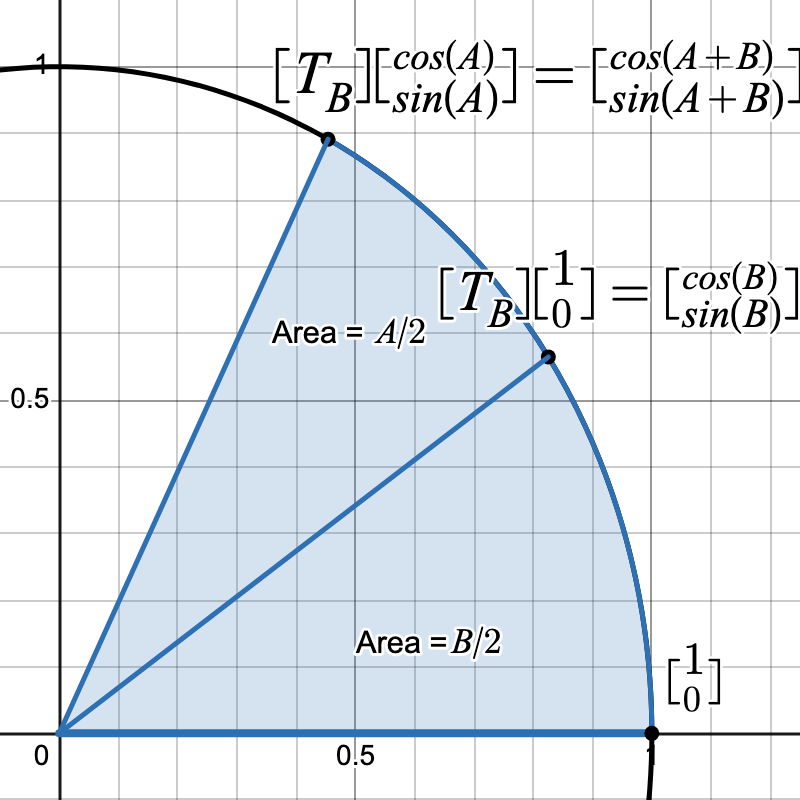}
    \caption{
    Using matrix $T_B = 
    \begin{bsmallmatrix}
    \cos(B) &-\sin(B) \\
    \sin(B) & \cos(B)
    \end{bsmallmatrix}$
    to add two areas $\frac{A}{2} + \frac{B}{2}$.   
    }
    \label{fig:Circ_sum_area}
\end{figure}

On the left-side diagram in Figure~\ref{fig:Circ_sum_area} we see a curvilinear triangle of area $A/2$ has been rotated counterclockwise to a similar region with the same area. The matrix for the linear map of rotation by an angle $B$ is determined by how it rotates the two basis elements 
$\begin{bsmallmatrix} 1 \\ 0 \end{bsmallmatrix}$ and 
$\begin{bsmallmatrix} 0 \\ 1 \end{bsmallmatrix}$.
The first is rotated to the point $\begin{bsmallmatrix}\cos(B) \\ \sin(B) \end{bsmallmatrix}$ and the second is rotated to the perpendicular point
$\begin{bsmallmatrix} -\sin(B) \\ \cos(B) \end{bsmallmatrix}$. Thus, the matrix
is given by the formula
\[
T_B=
\left[
\begin{array}{rr}
\cos(B) &-\sin(B) \\
\sin(B) & \cos(B)
\end{array}
\right].
\]
The determinant of this matrix is one, since $\det T_B =\cos^2(B) + \sin^2(B) =1 $, and so the transformation preserves area.

In the right-side diagram in Figure~\ref{fig:Circ_sum_area}, we have replaced the lower curvilinear triangle with one of area $B/2$ and we see its top edge aligns exactly with the edge of the rotated curvilinear triangle of area $A/2$, since 
$\begin{bsmallmatrix}
    \cos(B) \\ \sin(B)
\end{bsmallmatrix} =
T_B\begin{bsmallmatrix}
    1\\0 
\end{bsmallmatrix}$. The total area of these shaded regions is $(A+B)/2$ and by examining the top-most vertex on the region we conclude 
$\begin{bsmallmatrix} \cos(A+B) \\ \sin(A+B) \end{bsmallmatrix} =
T_B \begin{bsmallmatrix} \cos(A) \\ \sin(A) \end{bsmallmatrix}$.
Writing this out in matrix-vector form, we have
\[
\left[
\begin{array}{r}
\cos(A+B) \\
\sin(A+B) 
\end{array}
\right]
=T_B
\left[
\begin{array}{r}
\cos(A) \\
\sin(A)
\end{array}
\right]
=
\left[
\begin{array}{rr}
\cos(B) &-\sin(B) \\
\sin(B) & \cos(B)
\end{array}
\right]
\left[
\begin{array}{r}
\cos(A) \\
\sin(A)
\end{array}
\right].
\]
Multiplying out the matrix yields
\[
\left[
\begin{array}{r}
\cos(A+B) \\
\sin(A+B) 
\end{array}
\right]
=
\left[
\begin{array}{r}
\cos(A)\cos(B) -\sin(A)\cos(B) \\
\sin(A)\cos(B) + \cos(A)\sin(B)
\end{array}
\right]
.
\]
In equation form, this gives the summation rules for the trig functions, with 
\begin{align*}
    \cos(A+B) &= \cos(A)\cos(B) -\sin(A)\sin(B) \\
    \sin(A+B) &= \sin(A)\cos(B) + \cos(A)\sin(B)
\end{align*}
Setting $A = B$ yields the two double angle formulas
\begin{eqnarray*}
    \cos(2A) &= \cos^2(A)-\sin^2(A), \\
    \sin(2A) &= 2\sin(A)\cos(A).
\end{eqnarray*}

\subsection{Finding the derivatives at zero}

\begin{figure}
    \centering
    \includegraphics[width=0.7\linewidth]{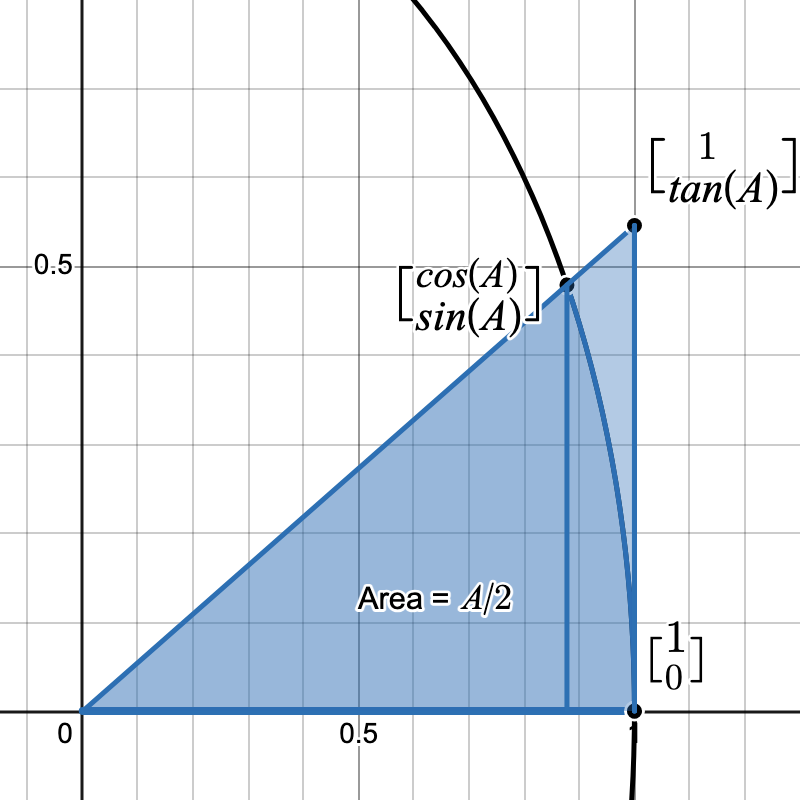}
    \caption{Regions for computing the derivative of $\sin(A)$ at $A=0$.}
    \label{fig:Circ_deriv}
\end{figure}

To find the derivatives, let's first consider the sine function at $A=0$. It is useful to guess the correct answer by estimation, before showing the result explicitly. Referring to Figure~\ref{fig:Circ_deriv}, we see the curvilinear section contains a smaller right triangle of height $\sin(A)$ and base $\cos(A)$ so its area is $\frac{1}{2}\sin(A)\cos(A)$. For small values of $A$, the cosine will be close to one and the area of the triangle is close to the area of the curvilinear section. Thus we have an approximation
\[ \frac{1}{2}\sin(A) \approx \frac{A}{2}, \]
more informatively
\[ \sin(A) \approx \sin(0) + A\cdot 1, \mbox{ for A small}.\]
So this looks like our linear approximation for $\sin(A)$ near $A=0$, showing we expect the derivative of sine to be equal to 1 here. 

To show this is correct linear approximation, we define $E(A) = \frac{\sin(A) - A}{|A|}$ for $A>0$ and $E(0) = 0$ giving the formula
\[ \sin(A) = \sin(0) + A\cdot 1 + |A|E(A).\]
We now verify the $E(A)$ is an error function. That is, it goes to zero as $A$ goes to zero. 

Referring again to Figure~\ref{fig:Circ_deriv}, the smaller triangle with height $\sin(A)$ is contained in the curvilinear region of area $A/2$ so we have the inequality of areas as
\[ \frac{1}{2}\sin(A)\cos(A) \leq \frac{A}{2}.\]
Using the double angle formulas from the previous section, we have $\frac{1}{2}\sin(2A) \leq A$ or equivalently $\sin(2A) \leq 2A$. Replacing $2A$ with $A$ gives the inequality $\sin(A) \leq A$ and thus
\[ \sin(A) - A \leq 0. \]
Next, by similar triangles, the larger triangle has height $\tan(A) = \frac{\sin(A)}{\cos(A)}$ and contains the curvilinear region of area $A/2$. So we have the area inequality $\frac{A}{2}\leq \frac{1}{2}\frac{\sin(A)}{\cos(A)}$ which we can rewrite as
\[A\cos(A)-A \leq \sin(A)- A. \]
Combining the two inequalities gives
\[A\cos(A)-A \leq \sin(A)- A \leq 0 \]
and dividing by $A>0$ gives
\[\cos(A)-1 \leq \frac{\sin(A)-A}{|A|} = E(A) \leq 0, \]
where the middle term is our error function $E(A).$

Now we see that as $A$ tends to zero, the difference $\cos(A)-1$ tends to $1-1 = 0$ by continuity of the cosine function, and so the bounding inequalities for $E(A)$ shows that it too tends to zero.

Thus $E(A)$ as defined is a proper error function, the approximation
\[ \sin(A) = \sin(0) + A\cdot 1 + |A|E(A)\]
is correct, and thus the derivation of $\sin(A)$ at $A=0$ is $\sin'(0) = 1.$

To find the derivative of $\cos(A)$ at $A=0$ we apply the chain rule to the formula $\cos(A) = \sqrt{1-\sin^2(A)}$ to find the derivative is $\cos'(0) = 0.$

\subsection{Finding the general derivative}

To find the derivatives at other points, we simply use the summation formulas from an earlier section. Treating $A=x$ as a variable and holding $B$ constant, we can differentiate the following formulas with respect to $x$
\begin{eqnarray*}
    \cos(x+B) &= \cos(x)\cos(B) -\sin(x)\sin(B) \\
    \sin(x+B) &= \sin(x)\cos(B) + \cos(x)\sin(B)
\end{eqnarray*}
to obtain
\begin{eqnarray*}
    \cos'(x+B) &= \cos'(x)\cos(B) -\sin'(x)\sin(B) \\
    \sin'(x+B) &= \sin'(x)\cos(B) + \cos'(x)\sin(B).
\end{eqnarray*}
Setting $x=0$ gives
\begin{eqnarray*}
    \cos'(B) &=& \cos'(0)\cos(B) -\sin'(0)\sin(B)  \\&=& 0\cdot\cos(B) -1\cdot\sin(B) \\ &=& -\sin(B),\\
    \sin'(B) &=& \sin'(0)\cos(B) + \cos'(0)\sin(B)  \\&=& 1\cdot\cos(B) + 0\cdot\sin(B)\\ &=& \cos(B).
\end{eqnarray*}

Thus the derivative of cosine is negative sine, and the derivative of sine is cosine. 

\section{The hyperbolic functions}

In this section, we define the hyperbolic functions $\cosh()$ and $\sinh()$ functions that arise in the geometry of the unit hyperbola in the $xy$-plane. We will consider the basic definitions for smaller parameter values, then extend them to arbitrarily large values, demonstrate the summation rules, and finally compute the derivatives. Our approach here is inspired by the geometric addition formulas derived in the short paper by Radcliffe~\cite{radcliffe_hyperbolic_addition}.

\subsection{Defining hyperbolic cosine and sine for small arguments}

\begin{figure}
    \centering
    \includegraphics[width=0.45\linewidth]{4_Hyper_defn_pos.png}
    \includegraphics[width=0.45\linewidth]{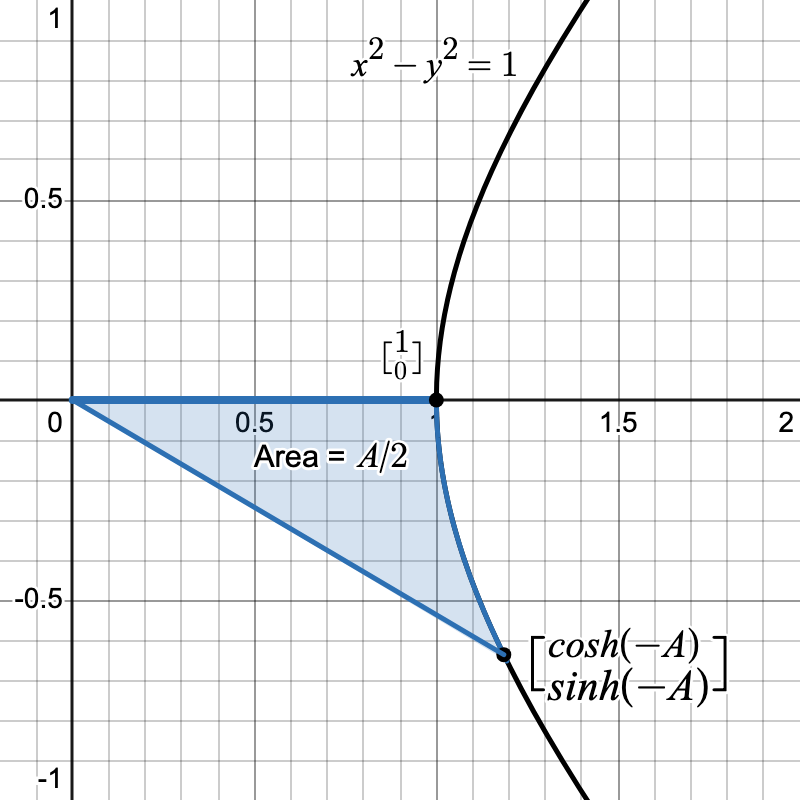}
    \caption{Given area $A/2$, hyperbolic functions $\cosh(A), \sinh(A)$ are defined.}
    \label{fig:Hyper_defn}
\end{figure}

We begin with the unit hyperbola $x^2 - y^2 =1$ on the plane, and focus our attention on the right branch, $x\geq 1$. Given a small positive number $A$, say less than 2 or so, we can form a curvilinear triangle of area $A/2$ with a vertex at the origin, a second vertex at the hyperbola's apex 
$\begin{bsmallmatrix}
    1\\0
\end{bsmallmatrix}$ 
and  a top vertex at a position 
$\begin{bsmallmatrix}
    x\\y
\end{bsmallmatrix}$
on the hyperbola, as shown in Figure~\ref{fig:Hyper_defn}. We  use this top vertex in the figure to define the values 
$\begin{bsmallmatrix}
    \cosh(A) = x\\\sinh(A) = y
\end{bsmallmatrix}$.

For the corresponding negative value $-A$, we drop a curvilinear triangle of area $A/2$ below the $x$-axis, as shown on the right side of Figure~\ref{fig:Hyper_defn}. This region has a lower vertex at a point 
$\begin{bsmallmatrix}
    x\\y
\end{bsmallmatrix}$
below the horizontal axis, which we use to define the hyperbolic cosine and sine as 
$\begin{bsmallmatrix}
    \cosh(-A) = x\\\sinh(-A) = y
\end{bsmallmatrix}$.

 From the symmetry in Figure~\ref{fig:Hyper_defn} we observe that including a negative sign in the argument does not change the value of the cosine, while it introduces a negative sign on the sine. Since the points are on the unit hyperbola, we observe that the difference of the squares of the hyperbolic cosine and sine is one. We have obtained three key hyperbolic relationships:
\begin{align*}
    \cosh(-A) &= \cosh(A) \\
    \sinh(-A) &= -\sinh(A) \\
    \cosh^2(A) - \sinh^2(A) &= 1 
\end{align*}

We can evaluate the hyperbolic functions at area zero by noting this identifies the apex of the hyperbola:

\begin{align*}
    \cosh(0) =&1,  &\sinh(0) =& 0   
\end{align*}

\subsection{Extending to large arguments}

For large values of $A>0$ it is not immediately apparent that a sufficiently large curvilinear region can be found with a large area $A/2$. 

\begin{figure}
    \centering
    \includegraphics[width=0.5\linewidth]{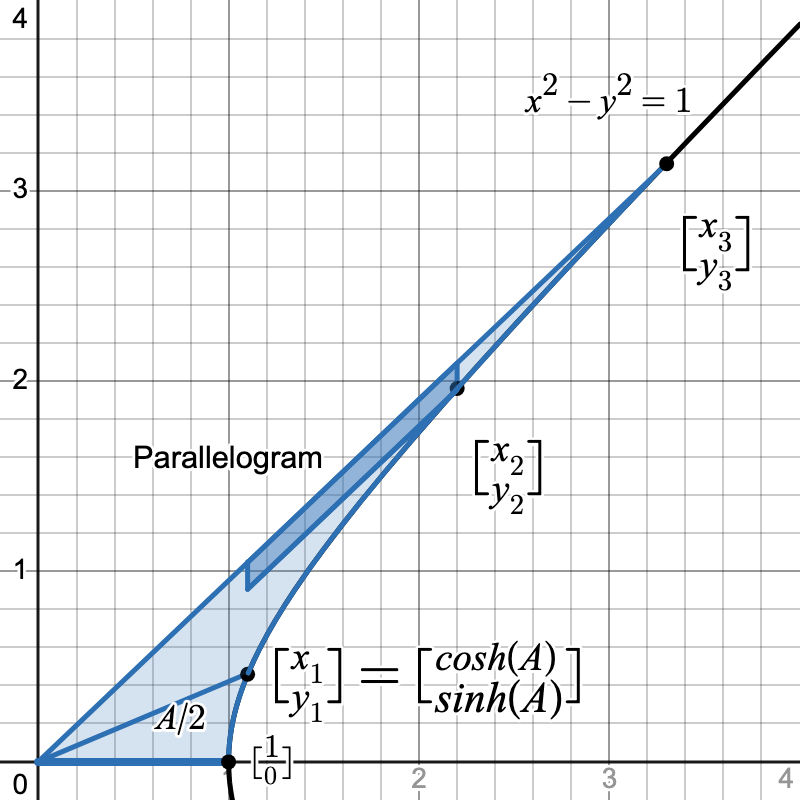}
    \caption{Increasing the area by at least a parallelogram.}
    \label{fig:Hyper_plus_pgram}
\end{figure}

To show that this is indeed the case, consider Figure~\ref{fig:Hyper_plus_pgram}. We say an area $A/2>0$ is {\it solvable} if we can find a point 
$\begin{bsmallmatrix}
    x\\y
\end{bsmallmatrix}$
on the unit hyperbola so the related curvilinear triangle has area $A/2$. This means $\cosh(A),\sinh(A)$ are defined. By the continuity of the area construction, if $A/2$ is solvable, so are all the smaller areas between 0 and $A/2$. Thus, it suffices to find arbitrarily large solvable areas to show all areas are solvable. 

Starting with some $x_1 = \cosh(A)$, we can choose points $x_3>x_2 >x_1$ and construct a certain parallelogram of area bigger than $1/9$, which always fits into the additional area for the region with vertex 
$\begin{bsmallmatrix}
    x_3\\y_3
\end{bsmallmatrix}$. 
This implies if $A/2$ is a solvable area, then $A/2 + 1/9$ is also a solvable area. By repeating this process, we can increase the solvable areas by multiples of 1/9 until we get to as large an argument as required.

The exact value $1/9$ is not critical. It arises from a convenient choice of $x$-coordinate values $x_1 = x, x_2 = 2x, x_3 = 3x$, as shown in Figure~\ref{fig:Hyper_plus_pgram}. The constructed parallelogram has vertical sides at $x=x_1$ and $x=x_2$ and two slanted sides: the top side along the line $y =(\frac{y_3}{x_3})x$ and the bottom side a parallel line passing through the point 
$\begin{bsmallmatrix}
    x_2\\y_2
\end{bsmallmatrix}$.
The area of the parallelogram is the base times the height, where here the base is the vertical thickness $b= (\frac{y_3}{x_3})x_2 - y_2$ and the height is the horizontal extent $h=x_2-x_1$. The area is thus
\begin{align*}
    b\cdot h &=  (\frac{x_2}{x_3}y_3 - y_2)\cdot (x_2-x_1) \\
    &= \left(\frac{2x}{3x}\sqrt{(3x)^2 - 1} - \sqrt{(2x)^2-1} \right)\cdot (2x-x) \\
    &= \left(2\sqrt{9x^2 - 1} - 3\sqrt{4x^2-1} \right)\frac{x}{3} \\
    &= \frac{(2\sqrt{9x^2 - 1})^2 - (3\sqrt{4x^2-1})^2}{2\sqrt{9x^2 - 1} + 3\sqrt{4x^2-1}}\frac{x}{3} \\
    &= \frac{4(9x^2-1) - 9(4x^2-1)}{2\sqrt{9x^2 - 1} + 3\sqrt{4x^2-1}}\frac{x}{3} 
    = \frac{5}{2\sqrt{9x^2 - 1} + 3\sqrt{4x^2-1}}\frac{x}{3} \\
    &> \frac{5}{2\sqrt{9x^2} + 3\sqrt{4x^2}}\frac{x}{3} =\frac{5}{12x}\frac{x}{3}
    = \frac{5}{36} > \frac{1}{9}.
\end{align*}
So the area is at least as big as 1/9. 

We conclude that all areas $A/2$ solvable, so the functions $\cosh(A)$, $\sinh(A)$ are defined by the geometry for all real values $A$.

\subsection{Summation rules}

In this section, we derive the summation rules for determining the values of $\cosh(A+B)$ and $\sinh(A+B)$. 

\begin{figure}
    \centering
    \includegraphics[width=0.45\linewidth]{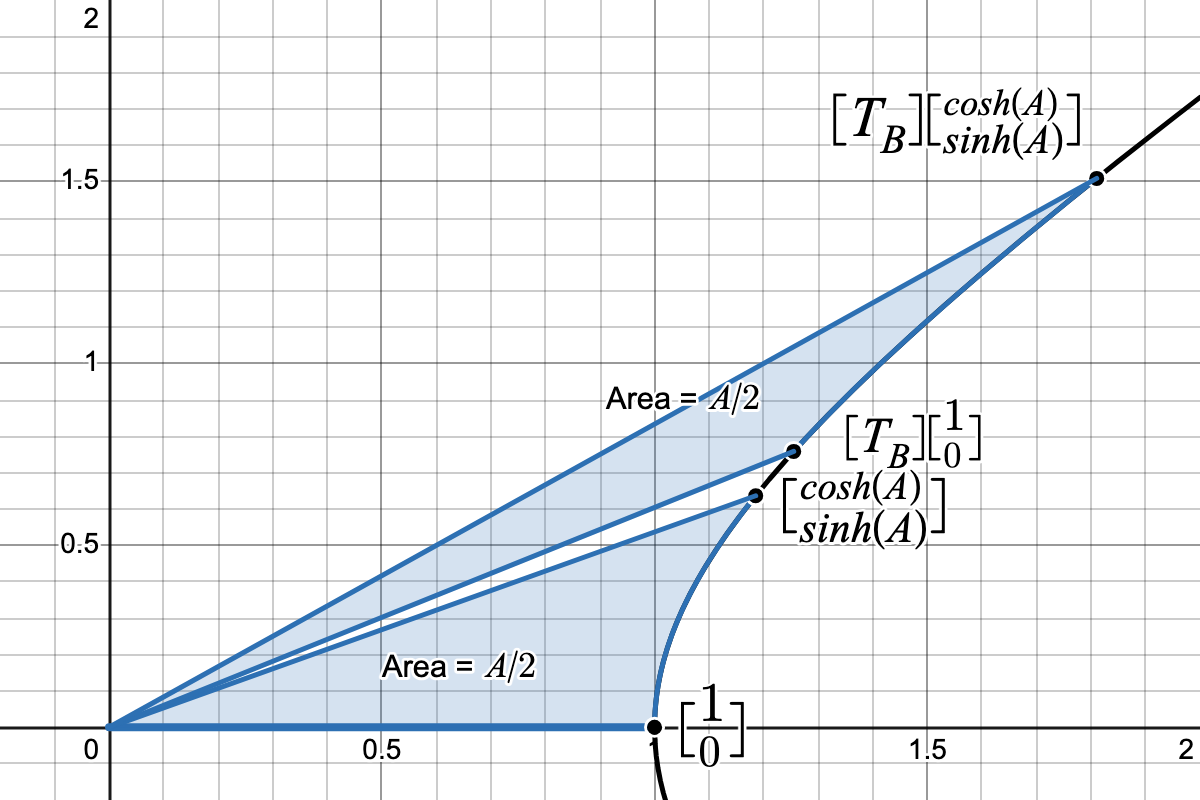}
    \includegraphics[width=0.45\linewidth]{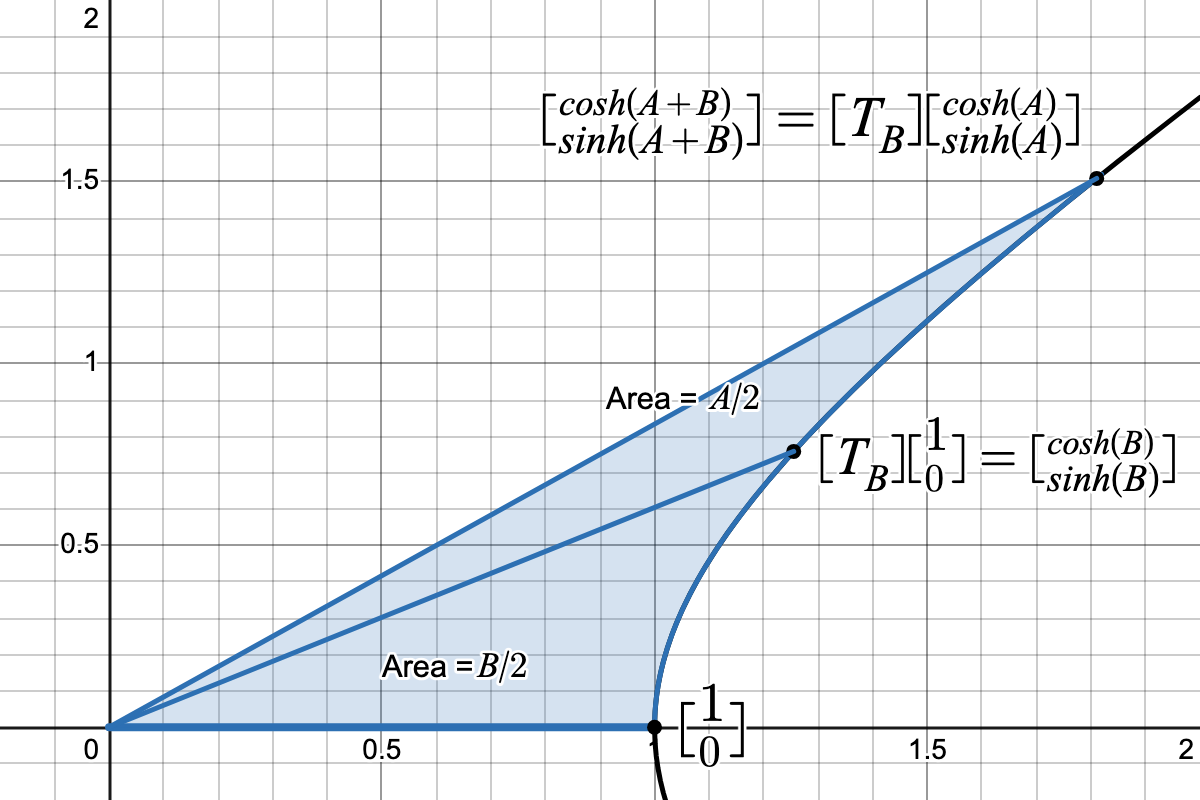}
    \caption{
    Using matrix $T_B = 
    \begin{bsmallmatrix}
    \cosh(B) &\sinh(B) \\
    \sinh(B) & \cosh(B)
    \end{bsmallmatrix}$
    to add two areas $\frac{A}{2} + \frac{B}{2}$.   
    }
    \label{fig:Hyper_sum_area}
\end{figure}

On the left-side diagram in Figure~\ref{fig:Hyper_sum_area} we see a curvilinear triangle of area $A/2$ has been rotated counterclockwise  and shifted onto a similar region with the same area. The matrix for this sheer transformation is determined by noting the point $\begin{bsmallmatrix} 1\\0\end{bsmallmatrix}$ at the apex must map to the point $\begin{bsmallmatrix}\cosh(B)\\ \sinh(B)\end{bsmallmatrix}$, while the corresponding point on the lower curve, $\begin{bsmallmatrix}\cosh(B)\\ -\sinh(B)\end{bsmallmatrix}$, must map to the apex $\begin{bsmallmatrix} 1\\0\end{bsmallmatrix}$. Solving, we find the matrix $T_B$ is given by the formula
\[
T_B=
\left[
\begin{array}{rr}
\cosh(B) &\sinh(B) \\
\sinh(B) & \cosh(B)
\end{array}
\right].
\]
The determinant of this matrix is one, since $\det T_B =\cosh^2(B) - \sinh^2(B) =1 $, and so the transformation preserves area. It is also easy to check this transformation maps the hyperbola onto itself. 

In the right-side diagram in Figure~\ref{fig:Hyper_sum_area}, we have replaced the lower curvilinear triangle with one of area $B/2$ and we see it aligns its edge exactly with the edge of the rotated curvilinear triangle of area $A/2$, since
$\begin{bsmallmatrix}\cosh(B)\\ \sinh(B)\end{bsmallmatrix} = 
T_B\begin{bsmallmatrix}1\\ 0
\end{bsmallmatrix}$. The total area of these shaded regions is $(A+B)/2$ and by examining the top-most vertex on the region we conclude $\begin{bsmallmatrix}\cosh(A+B)\\ \sinh(A+B)\end{bsmallmatrix} = T_B\begin{bsmallmatrix}\cosh(A)\\ \sinh(A)\end{bsmallmatrix}$. Writing this out in matrix-vector form, we have
\[
\left[
\begin{array}{r}
\cosh(A+B) \\
\sinh(A+B) 
\end{array}
\right]
=T_B
\left[
\begin{array}{r}
\cosh(A) \\
\sinh(A)
\end{array}
\right]
=
\left[
\begin{array}{rr}
\cosh(B) &\sinh(B) \\
\sinh(B) & \cosh(B)
\end{array}
\right]
\left[
\begin{array}{r}
\cosh(A) \\
\sinh(A)
\end{array}
\right].
\]
Multiplying out the matrix yields
\[
\left[
\begin{array}{r}
\cosh(A+B) \\
\sinh(A+B) 
\end{array}
\right]
=
\left[
\begin{array}{r}
\cosh(A)\cosh(B) +\sinh(A)\cos(B) \\
\sinh(A)\cosh(B) + \cosh(A)\sin(B)
\end{array}
\right]
.
\]
In equation form, this gives the summation rules for the trig functions, with 
\begin{align*}
    \cosh(A+B) &= \cosh(A)\cosh(B) +\sinh(A)\sinh(B) \\
    \sinh(A+B) &= \sinh(A)\cosh(B) + \cosh(A)\sinh(B)
\end{align*}
Setting $A = B$ yields the two double angle formulas
\begin{eqnarray*}
    \cosh(2A) &= \cosh^2(A)+\sinh^2(A), \\
    \sinh(2A) &= 2\sinh(A)\cosh(A).
\end{eqnarray*}

Note that the summation rules provide another explanation for why we can define $\cosh(),\sinh()$ for arbitrarily large arguments. For if we find a point $(x,y) = (\cosh(A),\sinh(A))$ on the unit hyperbola that gives a region of area $A/2$, then the point 
\[(x^2 + y^2,2xy) = (\cosh^2(A)+\sinh^2(A), 2\sinh(A)\cosh(A))\]
will give a region of area $A = 2(A/2)$. So we can always double the size of solvable area, and repeat over and over to get large enough areas to cover all values of the argument $A$. 

\subsection{Finding the derivatives at zero}

\begin{figure}
    \centering
    \includegraphics[width=0.7\linewidth]{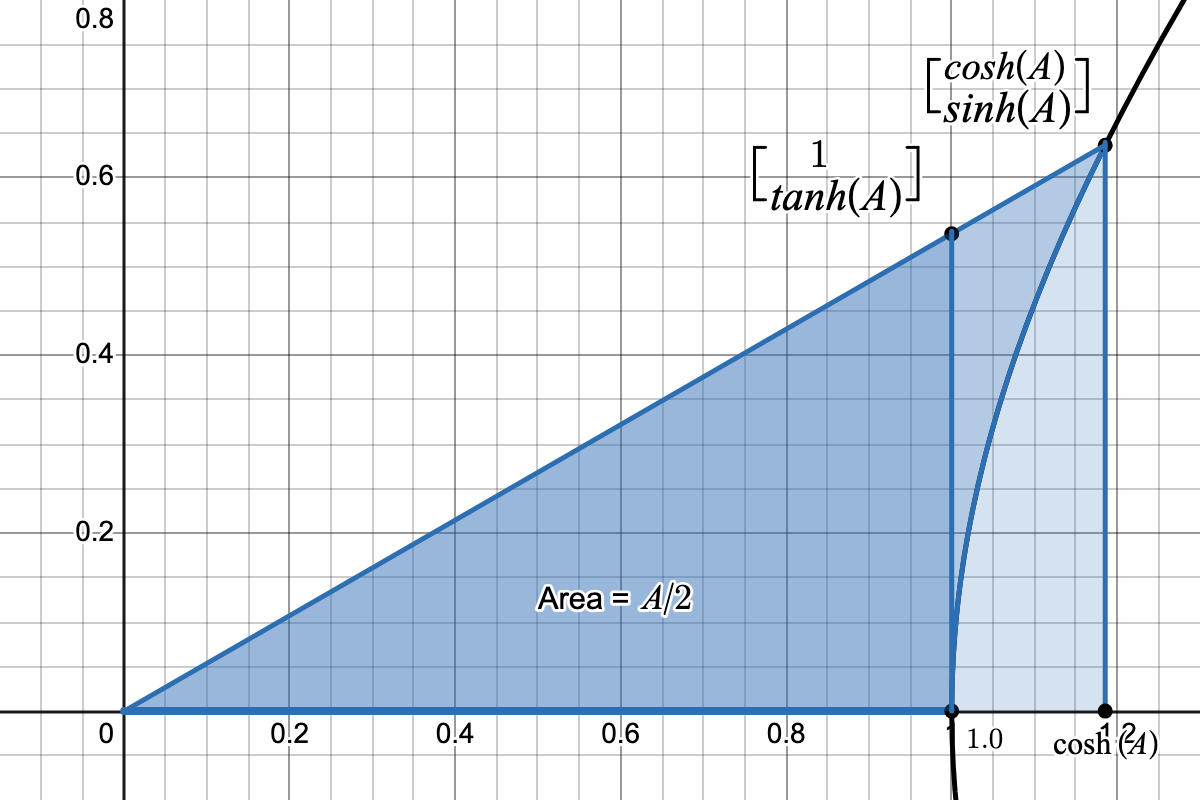}
    \caption{Regions for computing the derivative of $\sinh(A)$ at $A=0$.}
    \label{fig:Hyper_deriv}
\end{figure}

To find the derivatives, let's first consider the hyperbolic sine function at $A=0$. It is useful to guess the correct answer by estimation, before showing the result explicitly. Referring to Figure~\ref{fig:Hyper_deriv}, we see the curvilinear section is contained in a larger right triangle of height $\sinh(A)$ and base $\cosh(A)$, with  area $\frac{1}{2}\sinh(A)\cosh(A)$. For small values of $A$, the hyperbolic cosine will be close to one and the area of the triangle is close to the area of the curvilinear section. Thus we have an approximation
\[ \frac{1}{2}\sinh(A) \approx \frac{A}{2}, \]
more informatively
\[ \sinh(A) \approx \sinh(0) + A\cdot 1, \mbox{ for A small}.\]
So this looks like our linear approximation for $\sinh(A)$ near $A=0$, showing we expect the derivative of hyperbolic sine to be equal to 1 here. 

To show this is the correct linear approximation, we define $E(A) = \frac{\sinh(A) - A}{|A|}$ for $A>0$ and $E(0) = 0$ giving the formula
\[ \sinh(A) = \sinh(0) + A\cdot 1 + |A|E(A).\]
We now verify that $E(A)$ is an error function. That is, it goes to zero as $A$ goes to zero. 

Referring again to Figure~\ref{fig:Hyper_deriv}, the larger triangle with height $\sinh(A)$ contains the curvilinear region of area $A/2$ so we have the inequality of areas as
\[ \frac{1}{2}\sinh(A)\cosh(A) \geq \frac{A}{2}.\]
Using the double angle formulas from the previous section, we have $\frac{1}{2}\sinh(2A) \geq A$ or equivalently $\sinh(2A) \geq 2A$. Replacing $2A$ with $A$ gives the inequality $\sinh(A) \geq A$ and thus
\[ \sinh(A) - A \geq 0. \]
Next, by similar triangles, the smaller triangle has height $\tanh(A) = \frac{\sinh(A)}{\cosh(A)}$ and is contained inside the curvilinear region of area $A/2$. So we have the area inequality $\frac{A}{2}\geq \frac{1}{2}\frac{\sinh(A)}{\cosh(A)}$ which we can rewrite as
\[A\cosh(A)-A \geq \sinh(A)- A. \]
Combining the two inequalities gives
\[A\cosh(A)-A \geq \sinh(A)- A \geq 0 \]
and dividing by $A>0$ gives
\[\cosh(A)-1 \geq \frac{\sinh(A)-A}{|A|} = E(A) \leq 0, \]
where the middle term is our error function $E(A).$

Now we see that as $A$ tends to zero, the difference $\cosh(A)-1$ tends to $1-1 = 0$ by continuity of the hyperbolic cosine function, and so the bounding inequalities for $E(A)$ show that it too tends to zero.

Thus $E(A)$ as defined is a proper error function, the approximation
\[ \sinh(A) = \sinh(0) + A\cdot 1 + |A|E(A)\]
is correct, and thus the derivation of $\sinh(A)$ at $A=0$ is $\sinh'(0) = 1.$

To find the derivative of $\cosh(A)$ at $A=0$ we apply the chain rule to the formula $\cosh(A) = \sqrt{1+\sinh^2(A)}$ to find the derivative is $\cosh'(0) = 0.$

\subsection{Finding the general derivative}

To find the derivatives at other points, we simply use the summation formulas from an earlier section. Treating $A=x$ as a variable and holding $B$ constant, we can differentiate the following formulas with respect to $x$
\begin{eqnarray*}
    \cosh(x+B) &= \cosh(x)\cosh(B) +\sinh(x)\sinh(B) \\
    \sinh(x+B) &= \sinh(x)\cosh(B) + \cosh(x)\sinh(B)
\end{eqnarray*}
to obtain
\begin{eqnarray*}
    \cosh'(x+B) &= \cosh'(x)\cosh(B) +\sinh'(x)\sinh(B) \\
    \sinh'(x+B) &= \sinh'(x)\cosh(B) + \cosh'(x)\sinh(B).
\end{eqnarray*}
Setting $x=0$ gives
\begin{eqnarray*}
    \cosh'(B) &=& \cosh'(0)\cosh(B) +\sinh'(0)\sinh(B)  \\&=& 0\cdot\cosh(B) +1\cdot\sinh(B) \\ &=& \sinh(B),\\
    \sinh'(B) &=& \sinh'(0)\cosh(B) + \cosh'(0)\sinh(B)  \\&=& 1\cdot\cosh(B) + 0\cdot\sinh(B)\\ &=& \cosh(B).
\end{eqnarray*}

Thus the derivative of $\cosh()$ is $\sinh()$ and  the derivative of $\sinh()$ is $\cosh()$.

\section{The exponential function}

In this section, we define the exponential function $\exp()$ using the geometry of a rotated unit hyperbola. Of course, this function can be defined directly as the sum of the two hyperbolic functions $\cosh()$ and $\sinh()$, but we find this geometric approach instructive. We will consider the basic definitions for smaller parameter values, then extend them to arbitrarily large values, demonstrate the summation rules, and finally compute the derivatives.

\subsection{Defining exponential for small arguments}

\begin{figure}
    \centering
    \includegraphics[width=0.45\linewidth]{4_exp_defn_pos.png}
    \includegraphics[width=0.45\linewidth]{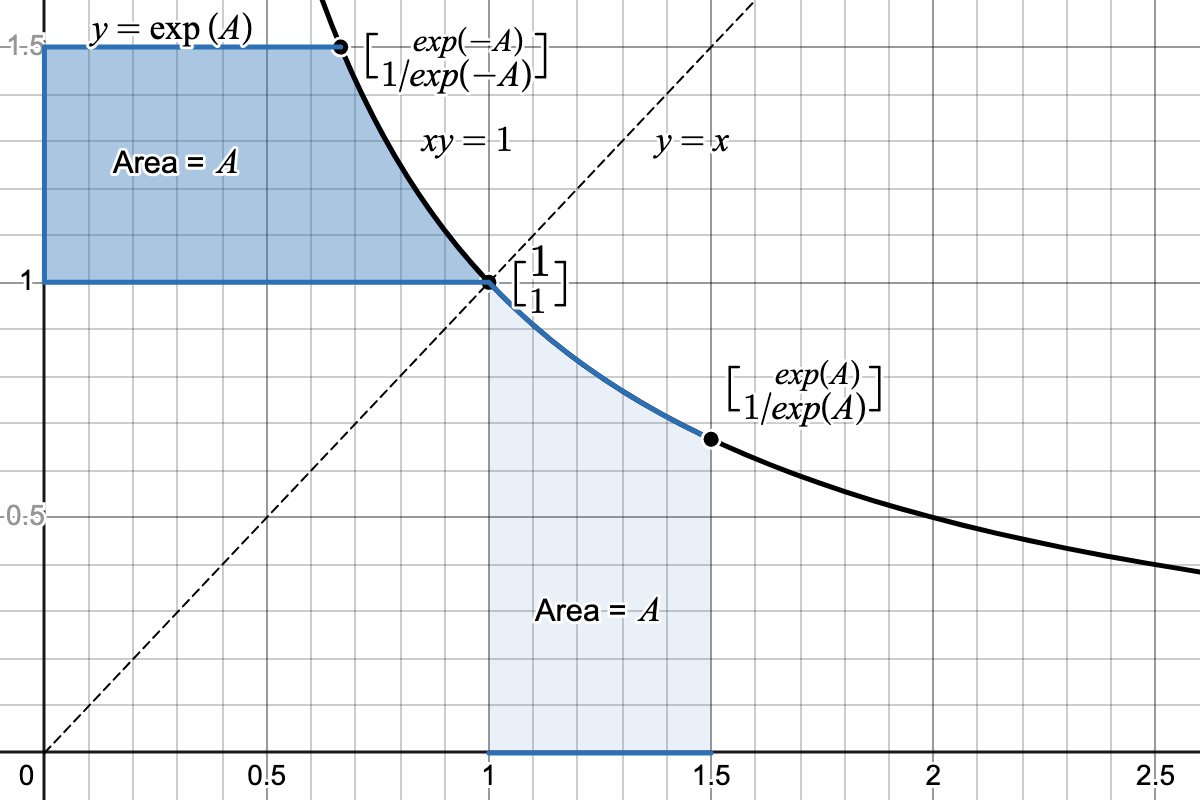}
    \caption{For area $A>0$, the values $\exp(A)$, $\exp(-A)$ are defined.}
    \label{fig:exp_defn}
\end{figure}

We begin with the (skew)\footnote{The term ``skew'' is not standard.} rectangular hyperbola $xy =1$ on the plane, and focus our attention the branch in the first quadrant, $x,y > 0$. Refer to the left-hand side of Figure~\ref{fig:exp_defn} for a geometric construction.  Given a small positive number $A$, say less than 2 or so, we construct a region of area $A$ that lies under the hyperbola and above the $x$-axis, and between the vertical lines $x=1$ and $x=x_0$. By continuity of the geometry, we can adjust the value of $x_0$ to get a region with exactly the area $A$. This point 
$\begin{bsmallmatrix}x_0\\ 1/x_0\end{bsmallmatrix}$
defines the exponential as $\exp(A) = x_0$. 

For the corresponding negative value $-A$, consider the right-hand side of Figure~\ref{fig:exp_defn}. We build a region between the $y$-axis and the hyperbola, and between the values $y=1$ and $y=y_0$. Again, we can adjust the value of $y_0$ to get exactly an area of $A$ for this region. The point 
$\begin{bsmallmatrix}1/y_0\\ y_0 \end{bsmallmatrix}$ 
defines the exponential for the negative value as $\exp(-A) = 1/y_0$, which is again the $x$-component of this point. 

By reflecting along the line $y=x$, it is clear the two regions in the figure have the same area, and we can conclude that
\[ \exp(-A) = \frac{1}{\exp(A)}. \]

\begin{figure}
    \centering
    \includegraphics[width=0.7\linewidth]{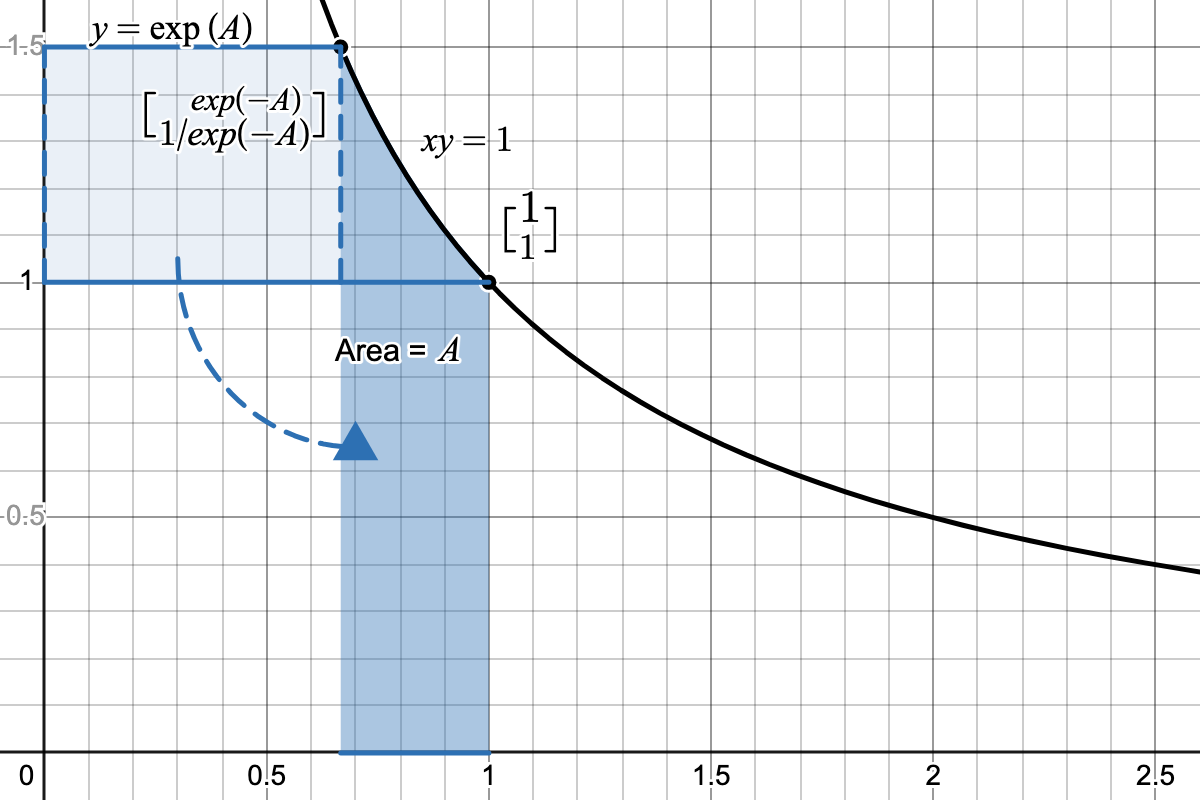}
    \caption{Alternate definition of $\exp(-A)$.}
    \label{fig:exp_defn_neg_neg}
\end{figure}

It is interesting to notice that the values of $\exp(-A)$ can also be defined as in Figure~\ref{fig:exp_defn_neg_neg}, where we consider a region under the hyperbola, bounded horizontally by $x_0 < x < 1$. A quick calculation shows the rectangle on the top left in the figure has the same area as the rectangle below the line $y=1$, so the two curvilinear regions and the same area. We will not use this fact in the following. 

We can evaluate the exponential function at zero as this corresponds to a vertical line at $x=1$ of zero area. Thus:
\[ \exp(0) = 1.\]

\subsection{Extending to large arguments}

For large values of $A>0$ it is not immediately apparent that a big enough  region can be found with large area $A$. 

\begin{figure}
    \centering
    \includegraphics[width=0.7\linewidth]{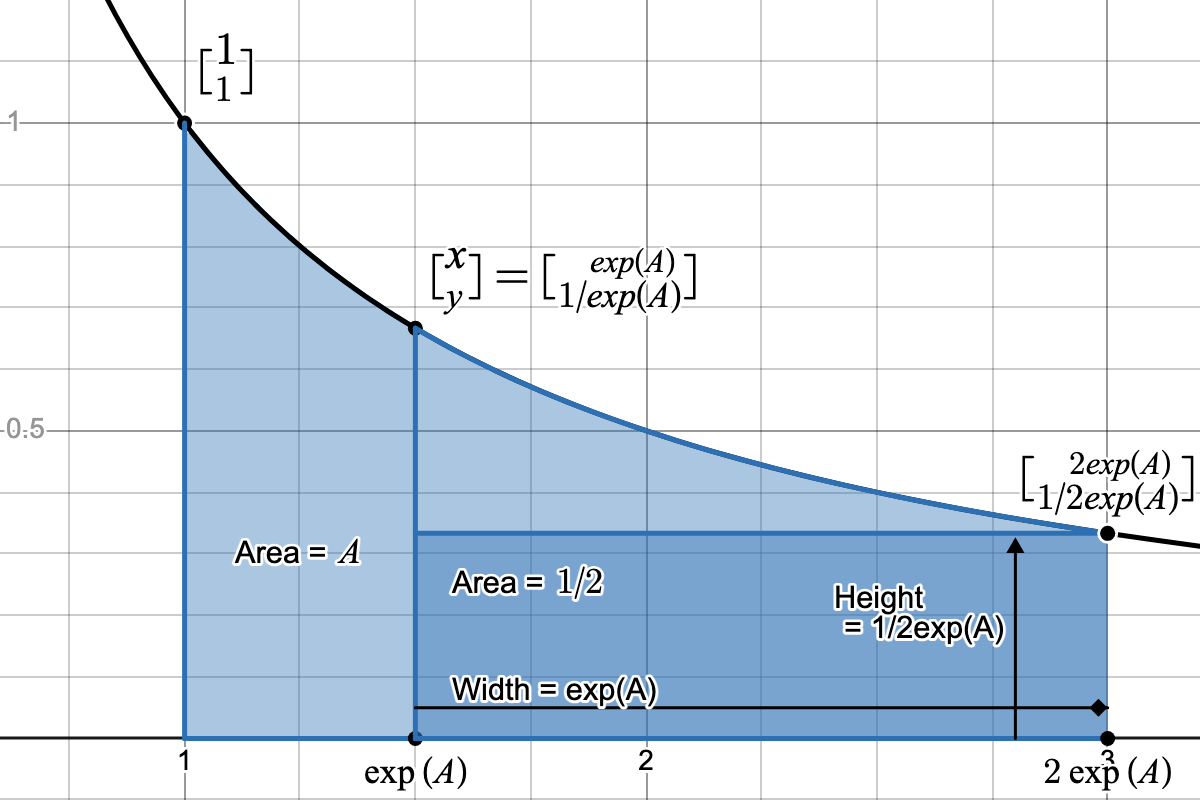}
    \caption{Increasing the area by at least a rectangle of area 1/2.}
    \label{fig:exp_plus_rect}
\end{figure}

To show this is indeed the case, consider Figure~\ref{fig:exp_plus_rect}. We say an area $A>0$ is {\it solvable} if we can find a point 
$\begin{bsmallmatrix}x\\ y  \end{bsmallmatrix}$ 
on the hyperbola so the related region has area $A$. (This means $\exp(A)$ is defined.) By continuity of the area construction, if some area $A>0$ is solvable, so are all the smaller areas. So it is enough to find arbitrarily large solvable areas to show all areas are solvable. 

Starting with some $x = \exp(A)$, we can choose points 
$\begin{bsmallmatrix}2x\\ 1/(2x)  \end{bsmallmatrix}$ 
on the hyperbola so that the area under the curve, from $x$ to $2x$ contains a rectangle of width $x$ and height 1$/(2x)$. This rectangle has area 1/2.  This implies if $A$ is a solvable area, then $A + 1/2$ is also a solvable area, and so is everything in between. By repeating this process, we can increase the solvable areas by multiples of 1/2 until we get to as large an argument as required.

We conclude that all areas $A>0$ solvable, so the function $\exp(A)$ defined by the geometry for all real values $A$.

\subsection{Summation rules}

In this section, we derive the summation rule for the exponential function. 

\begin{figure}
    \centering
    \includegraphics[width=0.45\linewidth]{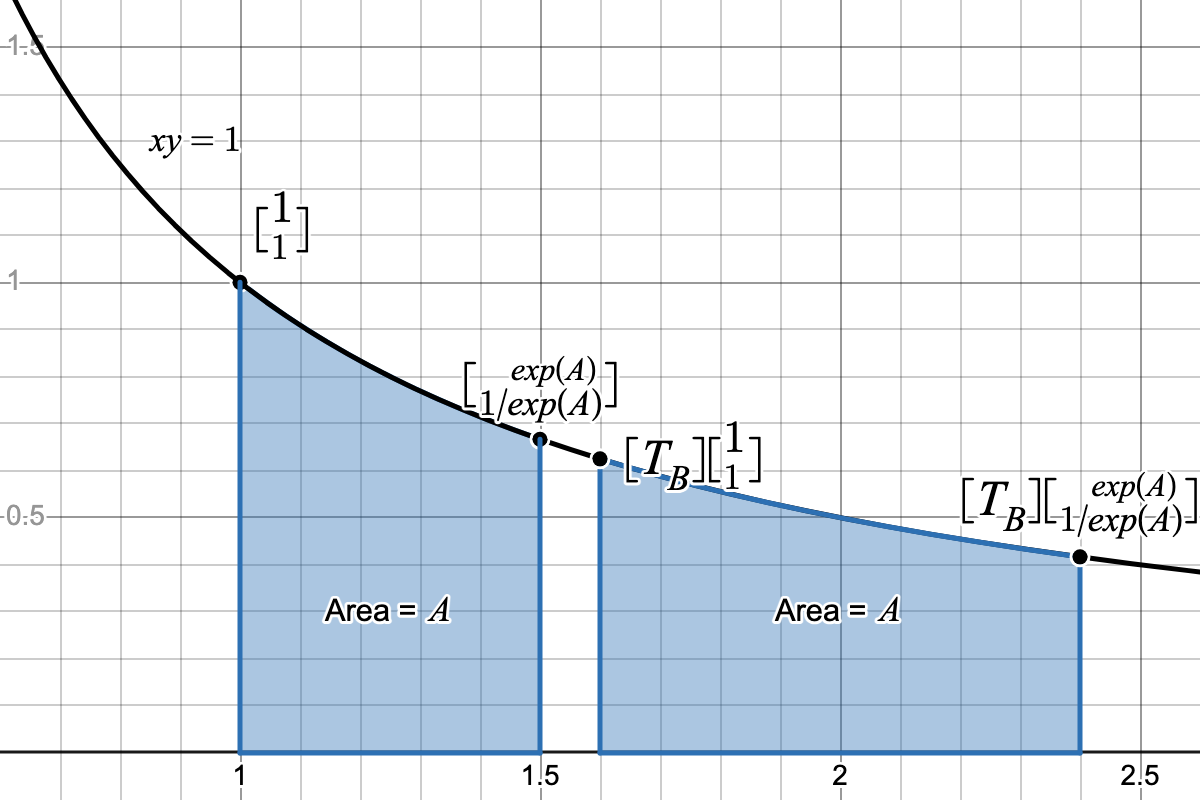}
    \includegraphics[width=0.45\linewidth]{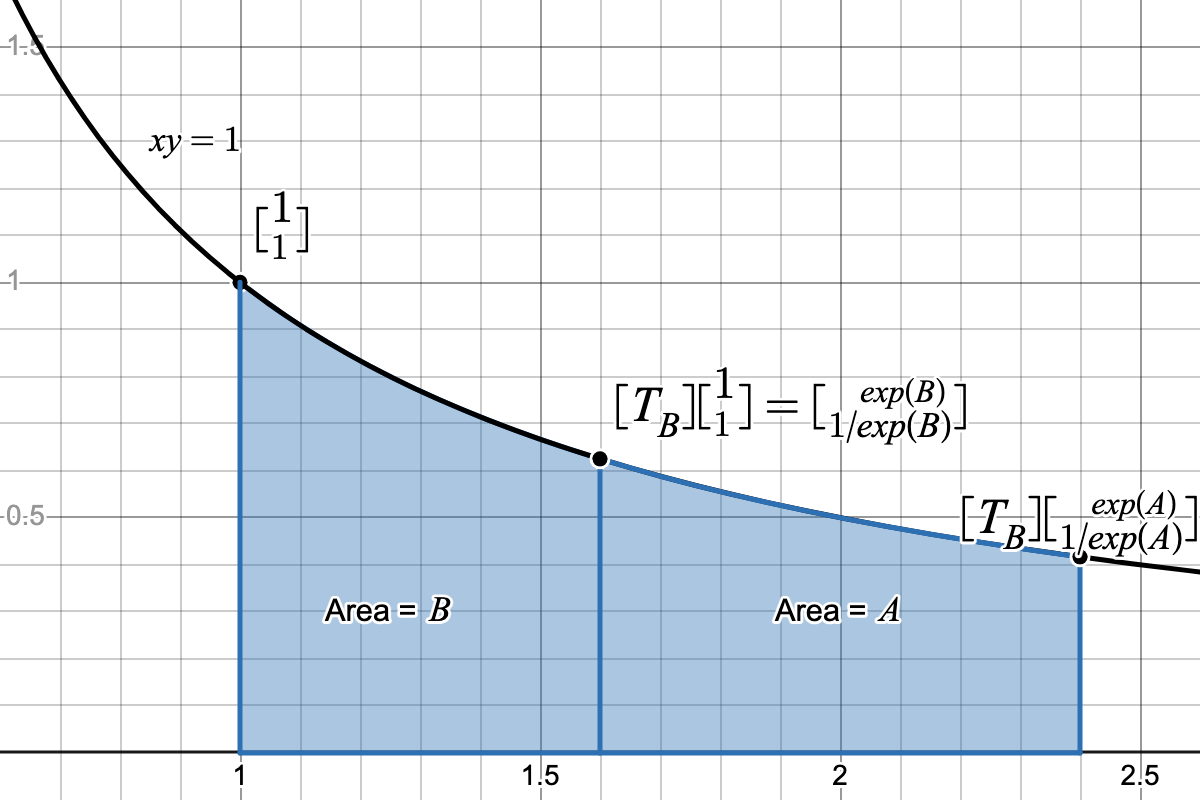}
    \caption{
    Using matrix $T_B = 
    \begin{bsmallmatrix}
    \exp(B) &0 \\
    0 & \exp(-B)
    \end{bsmallmatrix}$
    to add two areas $A+B$.   
    }
    \label{fig:exp_sum_area}
\end{figure}

On the left-side diagram in Figure~\ref{fig:exp_sum_area} we see a curvilinear region of area $A$ has been shifted right and compressed vertically into a similar region with the same area. The matrix for this linear transformation is determined by the mapping of the  vectors 
$\begin{bsmallmatrix}1\\ 0  \end{bsmallmatrix}$ 
and 
$\begin{bsmallmatrix}1\\ 1  \end{bsmallmatrix}$,
and is given by the formula
\[
T_B=
\left[
\begin{array}{cc}
\exp(B) &0 \\
0 & \exp(-B)
\end{array}
\right].
\]
The determinant of this matrix is one, since $\det T_B =\exp(B)\exp(-B) = \exp(B)/\exp(B) =1 $, and so the transformation preserves area. It is also easy to check this transformation maps the hyperbola $xy=1$ onto itself. 

In the right-side diagram in Figure~\ref{fig:exp_sum_area}, we have replaced the left region of area $A$ with a region of area $B$ and we see it aligns its edge exactly with the edge of the shifted region of area $A$,
since
$\begin{bsmallmatrix}\exp(B)\\ \exp(-B)\end{bsmallmatrix} = 
T_B\begin{bsmallmatrix}1\\ 1
\end{bsmallmatrix}$. 
The total area of these shaded regions is $A+B$ and by examining the right-most vertex on the regions we conclude 
$\begin{bsmallmatrix}\exp(A+B)\\ \exp(-A-B)\end{bsmallmatrix} = T_B\begin{bsmallmatrix}\exp(A)\\ \exp(-A)\end{bsmallmatrix}$.
Writing this out in matrix-vector form, we have
\[
\left[
\begin{array}{c}
\exp(A+B) \\
\exp(-A-B) 
\end{array}
\right]
=T_B
\left[
\begin{array}{c}
\exp(A) \\
\exp(-A)
\end{array}
\right]
=
\left[
\begin{array}{cc}
\exp(B) &0 \\
0 & \exp(-B)
\end{array}
\right]
\left[
\begin{array}{c}
\exp(A) \\
\exp(-A)
\end{array}
\right].
\]
Multiplying out the matrix yields
\[
\left[
\begin{array}{c}
\exp(A+B) \\
\exp(-A-B) 
\end{array}
\right]
=
\left[
\begin{array}{c}
\exp(A)\exp(B)  \\
\exp(-A)\exp(-B) 
\end{array}
\right]
.
\]
In equation form, this gives the addition rule for the exponential function, with 
\[ \exp(A+B) = \exp(A)\exp(B). \]

\subsection{Finding the derivative at zero}

\begin{figure}
    \centering
    \includegraphics[width=0.7\linewidth]{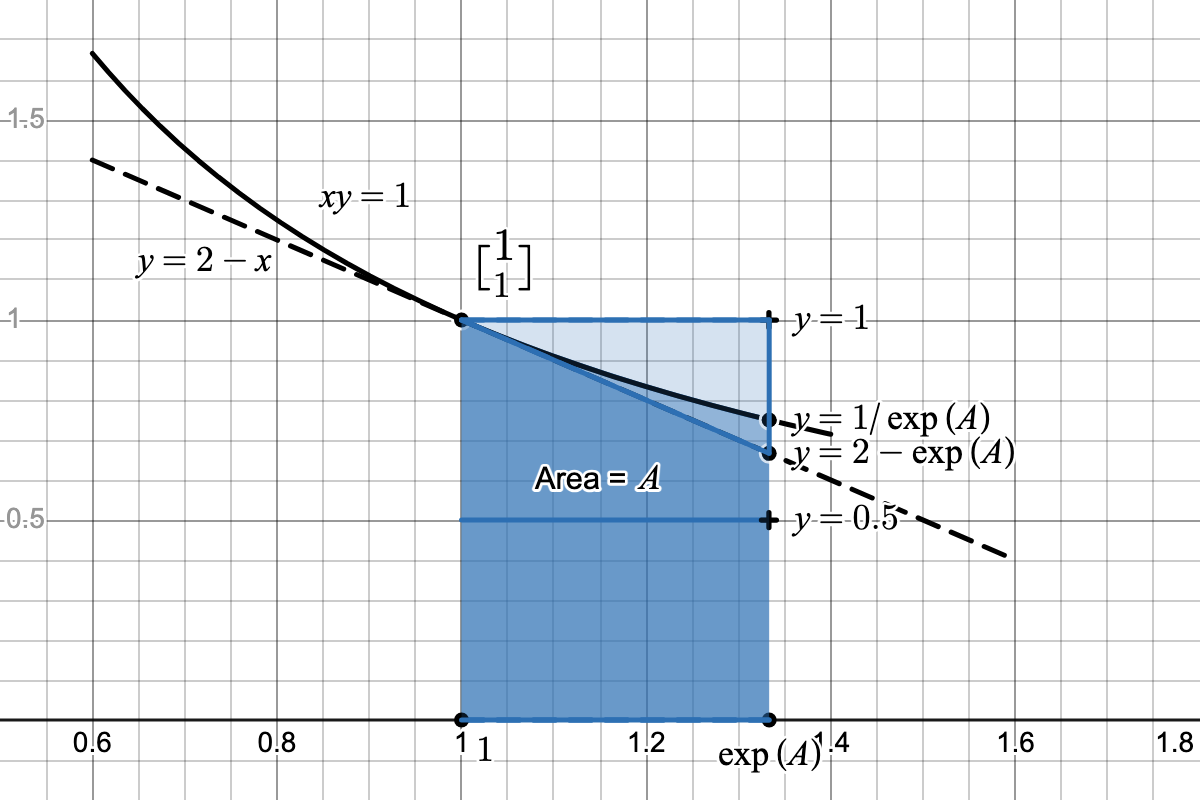}
    \caption{Regions for computing the derivative of $\exp(A)$ at $A=0$.}
    \label{fig:exp_deriv}
\end{figure}

To find the derivative, let's first consider the behaviour of the exponential function at $A=0$. It is useful to guess the correct answer by estimation, before showing the result explicitly. Referring to Figure~\ref{fig:exp_deriv}, we see the area of curvilinear section under the curve is close in area to the larger rectangle of height one and base $\exp(A)-1$. The area of this rectangle is $\exp(A) -1$, which is close to $A$, so we have an approximation
\[ \exp(A) - 1 \approx A. \]
More informatively
\[ \exp(A) \approx \exp(0) + A\cdot 1, \mbox{ for A small}.\]
So this looks like our linear approximation for $\exp(A)$ near $A=0$, showing we expect the derivative of the exponential to be equal to 1 there. 

To show this is the correct linear approximation, define $E(A) = \frac{\exp(A) - 1 - A}{|A|}$ for $A>0$, and set $E(0) = 0,$ giving the formula
\[ \exp(A) = \exp(0) + A\cdot 1 + |A|E(A).\]
We now verify that $E(A)$ is an error function. That is, it goes to zero as $A$ goes to zero. 

We note that the function $E(A) = \frac{\exp(A) - 1 - A}{|A|}$ is the ratio of two areas. The numerator $\exp(A) - 1 - A$ represents the area of the set difference between the large rectangular region shown in Figure~\ref{fig:exp_deriv}, whose area is $\exp(A)-1$, and the region of area $A$. This set difference is contained in the triangle formed from the tangent line $y=2-x$, with vertices at 
$\begin{bsmallmatrix} 1\\1 \end{bsmallmatrix}$,
$\begin{bsmallmatrix} \exp(A)\\2-\exp(A) \end{bsmallmatrix}$,
and
$\begin{bsmallmatrix} \exp(A)\\1 \end{bsmallmatrix}$.
The area of this set difference is less than the area of the triangle, which has a base and height both of length $\exp(A)-1$. So the inequality of areas gives
\[ 0\leq \exp(A)-1-A \leq \frac{1}{2}(\exp(A)-1)^2.\] 
Similarly, the denominator of the function $E(A)$ is the area $A$, which contains the smaller rectangle of height 1/2 and base $\exp(A)-1$, as shown in Figure~\ref{fig:exp_deriv}. This gives an inequality of areas,
\[ A \geq \frac{1}{2}(\exp(A)-1) > 0. \]
Combining these last two inequalities gives
\[ 0\leq E(A) = \frac{\exp(A) - 1 - A}{|A|} \leq \frac{\frac{1}{2}(\exp(A)-1)^2}{\frac{1}{2}(\exp(A)-1)} = \exp(A) - 1. 
\]

Now, as $A$ tends to zero, the difference $\exp(A)-1$ tends to zero, by continuity, so indeed $E(A)$ tends to zero as well. Thus, the function $E(A)$ is indeed an error function, as required.
The approximation
\[ \exp(A) = \exp(0) + A\cdot 1 + |A|E(A)\]
is correct, and thus the derivative of $\exp(A)$ at $A=0$ is $\exp'(0) = 1.$

\subsection{Finding the general derivative}

To find the derivative at other points, we simply use the summation formula for the exponential. Treating $A=x$ as a variable and holding $B$ constant, we can differentiate the following formula with respect to $x$:
\[ 
    \exp(x+B) = \exp(x)\exp(B)) \]
to obtain:
\[
    \exp'(x+B) = \exp'(x)\exp(B).  
\]
Setting $x=0$ gives
\[
    \exp'(B) = \exp'(0)\exp(B) = 1\cdot \exp(B).  
\]

Thus, the derivative of $\exp()$ is itself, $\exp()$.

\section{L'H\^{o}pital's Rule and $x\ln(x)$}

The well-known L'H\^{o}pital's Rule is traditionally expressed via limits, so it is important to verify that our no-limit approximation method still works here. We will show it does, even for the challenging cases of the rule. 

It is interesting to note that L'H\^{o}pital formulated his rule~\cite{bradley2015lhopital,lhopital1696analyse} long before Cauchy gave the formal definition of limits. It is not a big leap to express the rule in terms of our method of approximation. 

In the approximation format, L'H\^{o}pital's Rule states that if two differentiable functions $f,g$ are both zero at $x_0$, and $L$ is a fixed constant, it follows that:
\[ \mbox{ if } \frac{f'(x_0 + \epsilon)}{g'(x_0 +\epsilon)} \approx L,
\mbox{ then } \frac{f(x_0 + \epsilon)}{g(x_0 +\epsilon)} \approx L.\]
The use of approximations gives a quick rationale for why L'H\^{o}pital's rule should be true. If $f,g$ are differentiable functions with $f(x_0) = g(x_0) =0$ and the ratio $f'(x_0)/g'(x_0)$ is finite, then we have the approximations 
\begin{eqnarray*} 
f(x_0+\epsilon) &=& f(x_0) + \epsilon\cdot f'(x_0) + |\epsilon| E_1(\epsilon) = 0 + \epsilon\cdot f'(x_0) + |\epsilon| E_1(\epsilon), \\
g(x_0+\epsilon) &=& g(x_0) + \epsilon\cdot g'(x_0)\epsilon + |\epsilon| E_2(\epsilon) = 0 + \epsilon\cdot g'(x_0) + |\epsilon| E_2(\epsilon). 
\end{eqnarray*}
The ratio of functions then satisfies
\[ \frac{f(x_0 + \epsilon)}{g(x_0 + \epsilon)} = 
\frac{\epsilon\cdot f'(x_0) + |\epsilon| E_1(\epsilon)}{\epsilon\cdot g'(x_0) + |\epsilon| E_2(\epsilon)} 
= \frac{f'(x_0) + \sgn(\epsilon) E_1(\epsilon)}{g'(x_0) + \sgn(\epsilon) E_2(\epsilon)},\]
which, with a bit of algebra, leads to the approximation
\[ \frac{f(x_0 + \epsilon)}{g(x_0 + \epsilon)} = 
\frac{f'(x_0)}{g'(x_0)} +E(\epsilon),\]
for an error function $E(\epsilon)$ that can be derived exactly -- we skip the details.

In terms of traditional limits, this gives one form of L'H\^{o}pital's rule:
\[ \lim_{x\to x_0} \frac{f(x)}{g(x)} = \frac{f'(x_0)}{g'(x_0)}, \mbox{ if }f(x_0) = g(x_0) = 0.
\]

A challenging but standard problem in calculus is to find the limit, as $x$ tends to $0$, of the product $x\ln(x).$ In this case, the above form of L'H\^{o}pital's rule does not apply, because the function $\ln(x)$ and its derivatives are not defined at $x=0$ and its reciprocal $f(x) = 1/\ln(x)$, while it can be extended continuously to $f(0) = 0$, it is also not differentiable at $x=0$. 

We will need a more general form of L'H\^{o}pital's rule, in approximation form.

Here is a sequence of results to get to the general L'H\^{o}pital's rule.

\subsection{Critical Point Theorem} 
{\bf CPT Theorem: } If $f(x)$ is continuous on the interval $[a,b]$, differentiable on $(a,b)$ and $c$ is a maximizer or minimizer for $f$ in the interior of the interval, then $f'(c) = 0$.

This can be proved by considering the approximation at $x_0 = c$ with
\[ f(c + \epsilon) = f(c) + f'(c)\epsilon + |\epsilon|E(\epsilon).\]
If $f'(c) \neq 0$, taking $\epsilon\neq 0$ small enough that $|E(\epsilon)| <|f'(c)/2|$ shows the term $f'(c)\epsilon$ dominates, so $f(c \pm \epsilon)$ is either bigger or smaller than $f(c)$. Therefore, $f(c)$ is neither a max nor a min for the function. So, for $c$ to be a minimizer or maximizer, it must be that the derivative there is zero. 

\subsection{Mean Value Theorem} 
{\bf MVT Theorem: } If $f(x)$ is continuous on the interval $[a,b]$ and differentiable on $(a,b)$, then there is a point  $c$ in the interval $(a,b) $  with
\[ f'(c) = \frac{f(b) - f(a)}{b-a}.\] 

The proof of this follows immediately from Critical Point Theorem, for if we define a helper function
\[h(x) = f(x) - \frac{f(b)-f(a)}{b-a}x,\]
we see that $h(b) = h(a)$, hence, it is either constant, or it has an interior minimizer or maximizer at $x=c$. Its derivative  is zero, so $h'(c) = 0 = f'(c) - \frac{f(b) - f(a)}{b-a}$.

\subsection{Cauchy Mean Value Theorem} 
{\bf CMVT Theorem: } If $f(x),g(x)$ are continuous functions on the interval $[a,b]$ and differentiable on $(a,b)$, with $g'(x)$ never zero, then there is a point $c$ in the interval $(a,b) $  with
\[ \frac{f'(c)}{g'(c)} = \frac{f(b) - f(a)}{g(b)-(a)}.\] 

To demonstrate this, note from the previous Mean Value Theorem, there is a point $c'$ with 
\[ g'(c') = \frac{g(b)-g(a)}{b-a} \]
which is non-zero by assumption. Therefore, $g(b)-g(a)$ is not zero, so we can  define a helper function
\[h(x) = f(x) - \frac{f(b)-f(a)}{g(b)-g(a)}g(x).\]
Again, we check that  $h(b) = h(a)$, so it has an interior minimizer or maximizer $c$ where its derivative is zero. Thus $h'(c) = 0 = f'(c) - \frac{f(b) - f(a)}{g(b)-g(a)}g'(c)$. Dividing through by $g'(c)$ gives the desired equality.

\subsection{General L'H\^{o}pital's Rule}
{\bf L'H\^{o}pital's Theorem: } Suppose $f(x),g(x)$ are continuous and differentiable functions on an interval containing  on $x_0$, with $f(x_0) = g(x_0) = 0$ and $g'(x)\neq 0$ for $x$ near, but not equal to, $x_0$. If
\[ \frac{f'(x_0 + \epsilon)}{g'(x_0 +\epsilon)} = L + E(\epsilon) \]
for some number $L$ and all small $\epsilon\neq 0$, with $E(\epsilon)$ an error function, then
\[ \frac{f(x_0 + \epsilon)}{g(x_0 +\epsilon)} = L + E_1(\epsilon), \]
for the same number $L$ and possibly a different error function $E_1(\epsilon)$. 

The result holds when restricting to the one-sided cases $\epsilon>0$ or $\epsilon <0$. It also holds if $f(x),g(x)$ are differentiable for all points $x$ near $x_0$ but not necessarily at $x_0$. 

Finally, the result holds in the unbounded case: rather than $f(x_0) = g(x_0)=0$, we assume the functions $f,g$ are unbounded near $x_0$, with $\frac{1}{g(x_0+\epsilon)}$ an error function. That is, the  function $\frac{1}{g(x_0+\epsilon)}$ is approximately zero for small $\epsilon$.   

Informally, given the conditions of the theorem, we have the statement about approximations in $\epsilon$:
\[ \mbox{ If } \frac{f'(x_0 + \epsilon)}{g'(x_0 +\epsilon)} \approx L,
\mbox{ then } \frac{f(x_0 + \epsilon)}{g(x_0 +\epsilon)} \approx L.\]

To show this result in the case where $f(x_0) = g(x_0) = 0$, by assumption there is an interval $(x_0,x_1)$ where the derivative $g'(x)$ is non-zero. Therefore, $g(x) \neq 0$ in this interval, otherwise we would have $g(x) - g(x_0) =0$ and thus $g'(x') =0$ at some intermediate point $x'$, by the Mean Value theorem. Thus for small $\epsilon>0$ we can divide by the non-zero $g(x_0 + \epsilon)$ and then by the Cauchy Mean Value Theorem, we have
\[ \frac{f(x_0+\epsilon)}{g(x_0+\epsilon)} = \frac{f(x_0+\epsilon) - f(x_0)}{g(x_0+\epsilon) - g(x_0)}
=\frac{f'(x_0 + \epsilon')}{g'(x_0 + \epsilon')} = L + E(\epsilon'),\]
where $\epsilon'$ is some value $0<\epsilon' < \epsilon,$ with $c=x_0+\epsilon'$ the interior point of satisfying the Cauchy MVT. 

Now set $E_1(\epsilon) = E(\epsilon')$ and note that as $\epsilon$ decreases to zero, so does $\epsilon'$ and thus so does $E_1(\epsilon) = E(\epsilon').$ Thus $E_1(\epsilon)$ is an error function, proving the result. 

To demonstrate the case where $f,g$ are unbounded near $x_0$, define a remainder function $r(x)$ by the equation \[ \frac{f(x)}{g(x)} = L + r(x).\] We will show the function $E_1(\epsilon) = r(x_0 + \epsilon)$ is an error function, thus proving the result. 

To demonstrate this, we must show $E_1(\epsilon)$ can be made arbitrarily small when $\epsilon$ is small enough (Property~3 from Section~3). Fix a bound $B>0$. Since $E(\epsilon)$ in the statement of the theorem is an error function, there is some $\epsilon_0 > 0$ such that $0<\epsilon<\epsilon_0$ implies $|E(\epsilon)| < B/2$. By reducing $\epsilon_0$ if necessary, we may assume the derivative $g'(x)$ is non-zero for all $x$ in the interval $(x_0,x_0+\epsilon_0)$. Thus, $g(x)- g(x_0+\epsilon_0)\neq 0$ for all points $x$ in $(x_0,x_0+\epsilon_0)$. 

Now, for any $x=x_0+\epsilon$ in this interval $(x_0,x_0+\epsilon_0)$, the Cauchy MVT yields a point $c = x_0 + \epsilon'$ with
\[ \frac{f(x_0+\epsilon) - f(x_0+\epsilon_0)}{g(x_0+\epsilon) - g(x_0+\epsilon_0)}
=\frac{f'(x_0 + \epsilon')}{g'(x_0 + \epsilon')} = L + E(\epsilon'),\]
where $0<\epsilon' < \epsilon_0$.
Divide the numerator and denominator on the left hand side by $g(x_0 + \epsilon)$ to get 
\[  \frac{\frac{f(x_0+\epsilon)}{g(x_0 + \epsilon)} - \frac{f(x_0 +\epsilon_0)}{g(x_0+\epsilon)}}{1 - \frac{g(x_0 + \epsilon_0)}{g(x_0+\epsilon)}}
= L + E(\epsilon').\]
Moving the denominator to the right hand side and rearranging gives
\[  \frac{f(x_0+\epsilon)}{g(x_0 + \epsilon)}
= L +  E(\epsilon')  + \frac{1}{g(x_0+\epsilon)}[ f(x_0 + \epsilon_0) - (E(\epsilon')+L)g(x_0+\epsilon_0)].\]
We recognize the right hand side as $L+ E_1(\epsilon)$ by our definition of the remainder function. The second term $E(\epsilon')$ is bounded in absolute value by $B/2$, by the construction since $0 < \epsilon' <\epsilon_0$. The third term can be made smaller than $B/2$ since it is a product of the error function $\frac{1}{g(x_0+\epsilon)}$ times a bounded factor. (Note this second factor is bounded as $|E(\epsilon')|$ is bounded by $B/2$ and all the other parts are constant.) Thus, the sum giving $E_1(\epsilon)$, is bounded in absolute value by $B/2 + B/2 = B$ for $\epsilon$ sufficiently small, demonstrating Property~3 of an error function and concluding the proof. 

\section{The example of $x\ln(x)$}

Now we consider the case of the approximation of $x\ln(x)$, and we want show that 
\[x\ln(x) = 0 + E(x), \mbox{ for small } x>0. \]
That is, we need to show that  $x\ln(x)$ gets close to $L=0$ as $x$ gets close to zero. 

To apply L'H\^{o}pital's rule, we write
\[ x\ln(x) = \frac{\ln(x)}{\frac{1}{x}} = \frac{f(x)}{g(x)},\]
with $f(x) = \ln(x)$ and $g(x) = \frac{1}{x}$. The reciprocal of $g$, $\frac{1}{g(x)} = x$, is an error function for small $x$, as required by the statement of the Rule.   Then
\[ \frac{f'(x)}{g'(x)} = \frac{\frac{1}{x}}{-\frac{1}{x^2}} = -x = 0 + E(x), \mbox{ for small } x,\]
where here the error function is just $E(x) = -x.$ Thus, from L'H\^{o}pital's rule, we have
\[  x\ln(x) = \frac{\ln(x)}{\frac{1}{x}} = 0 + E_1(x),\]
for small $x$.

In other words, for $x$ close to 0, the function $x\ln(x)$ is close to zero. 

\section{Approximation with Taylor polynomials}

From L'H\^{o}pital's rule, we immediately obtain the Taylor polynomial approximation for n-times differentiable functions. We state the theorem as follows:

{\bf Taylor's Theorem:} Suppose the function $f(x)$ and its derivatives $f'(x)$, $f''(x)$, ..., $f^{(n-1)}(x)$ are defined in an interval around a point $x_0$, and the n-th derivative $f^{(n)}(x_0)$ exists at $x_0$. Then the n-th order polynomial 
\[p_n(\epsilon) = f(x_0) + f'(x_0)\epsilon + \frac{f^{(2)}(x_0)}{2!}\epsilon^2 +\frac{f^{(3)}(x_0)}{3!}\epsilon^3 +
\cdots + \frac{f^{(n)}(x_0)}{n!}\epsilon^n\]
approximates $f(x_0+\epsilon)$ well for small values of $\epsilon$, with
\[ f(x_0 + \epsilon) = p(\epsilon) + |\epsilon|^n E(\epsilon)\]
for some error function $E(\epsilon)$. 

Note the error function $E(\epsilon)$ is known as the Peano form of the remainder, up to sign. 

We produce a proof using the approximation approach, which follows from our definition of the n-th derivative and repeated application of L'H\^{o}pital's rule. To see this, first note that to define the higher order derivations, the functions $f, f' f'', \cdots f^{(n-1)}$ all need to be defined near the point $x_0$. However, the statement of the theorem requires the n-th derivative only to exist at the one point $x_0$. By our approximation definition for the derivative of the function $f^{(n-1)}$, there is an error function $E(\epsilon)$ satisfying the following equality:
\[ f^{(n-1)}(x_0+\epsilon) = f^{(n-1)}(x_0) + f^{(n)}(x_0)\epsilon + |\epsilon|E(\epsilon).\] 
Rearranging and dividing by $n!\epsilon$ gives
\[ \frac{f^{(n-1)}(x_0+\epsilon) - f^{(n-1)}(x_0)}{n!\epsilon} = \frac{f^{(n)}(x_0)}{n!} + E_1(\epsilon) \] 
where $E_1(\epsilon) =\frac{1}{n!}\sgn(\epsilon)E(\epsilon)$ is also an error function. 
This is now in the form to apply L'H\^{o}pital's rule, with $L = \frac{f^{(n)}(x_0)}{n!}$ playing the role of $L$ in the statement of the rule. 

We define two functions
\begin{eqnarray*}
    F(\epsilon) &=& f(x_0+\epsilon) - \sum_{k=0}^{n-1}\frac{f^k(x_0)}{k!} \epsilon^k, \\
    G(\epsilon) &=& \epsilon^n.
\end{eqnarray*}
Computing their derivatives, we find the functions $F,G$ and their first $n-1$ derivatives are all zero at $\epsilon=0$.  From L'H\^{o}pital's rule, we have
\[ \frac{F(\epsilon)}{G(\epsilon)} = \cdots = \frac{F^{(n-1)}(\epsilon)}{G^{(n-1)}(\epsilon)} + E_2(\epsilon).\]
Expanding the terms, we find
\[ \frac{f(x_0 + \epsilon) - \sum_{k=0}^{n-1} \frac{f^{(k)}(x_0)}{k!} \epsilon^k}{\epsilon^n}  = \cdots = \frac{f^{(n-1)}(x_0+\epsilon)-f^{(n-1)}(x_0)}{n!\epsilon} + E_2(\epsilon).\]
But the right hand side we have already seen above as the definition of the n-th derivative, so we have
\[ \frac{f(x_0 + \epsilon) - \sum_{k=0}^{n-1} \frac{f^{(k)}(x_0)}{k!} \epsilon^k}{\epsilon^n} = \frac{f^{(n)}(x_0)}{n!} + E_1(\epsilon) + E_2(\epsilon),
\] 
which can be rearranged as 
\[ f(x_0 + \epsilon) =  \left( \sum_{k=0}^{n-1} \frac{f^{(k)}(x_0)}{k!} \epsilon^k \right) + \frac{f^{(n)}(x_0)}{n!}\epsilon^n + |\epsilon|^n E_3(\epsilon).
\] 
Here we have the error function $E_3(\epsilon) = \sgn(\epsilon)(E_1(\epsilon) + E_2(\epsilon))$

In other words, for the polynomial $p(\epsilon)$ in the statement of the theorem, we have the approximation
\[ f(x_0 + \epsilon) =  p(\epsilon) + |\epsilon|^n E_3(\epsilon),
\]
which completes the proof. 

{\bf Remark:}
We note a good polynomial approximation for a function can exist even when it is not differentiable anywhere except at one point. For instance, with $w(x)$ the nowhere differentiable Weierstrauss function, the function
\[ f(x) = 1 + x + x^2 + x^3 + x^4 w(x)\] 
is differentiable only at the point $x_0 = 0$ yet it is very well approximated by the third order polynomial
\[ p(x) = 1 + x + x^2 + x^3.\]
This suggests we could define higher order derivatives for $f$ at an isolated point, even when the derivatives do not exist near that point, by examining the coefficients of the approximating polynomial. An interesting idea worth pursuing. 

\section{A peek at integral calculus}

Let us show how error functions are used to demonstrate the Fundamental Theorem of Calculus (FTC). That is, the integral of a continuous function can be computed from its antiderivative.

We gloss over some details for expediency here. For instance, for a non-negative function $f(x)$ we would define an integral $\int_{x_0}^{x_1} f(x)$ as the area under the curve $y=f(x)$, between the two endpoints $x_0$ and $x_1$. This can be made precise by noting for any piecewise constant functions $m(x), M(x)$ which bound $f(x)$ as $m(x) \leq f(x) \leq M(x)$ then by the containment of geometric regions, we have an inequality of areas with 
\[ \int_{x_0}^{x_1} m(x) dx \leq \int_{x_0}^{x_1} f(x) dx \leq \int_{x_0}^{x_1} M(x) dx,
\]
where the areas on the left and right sides of these inequalities are well-defined as the sum of areas of rectangles. We then use uniform continuity of $f(x)$ to show the supremum of the integrals on the LHS equals the infinum of the integrals on the RHS, which gives a well-defined value for the integral of $f(x)$. 

Now to demonstrate the Fundamental Theorem. Suppose $f(x)$ is continuous and we define the function $F(x)$ as the integral (that is, an area under a curve) as 
\[F(x_1) = \int_{x_0}^{x_1} f(x) \, dx.\]
Since $f$ is continuous at $x_1$, we can write $f(x_1 + t) = f(x_1) + E(t)$ for some error function $E(t)$, where $t$ is a small parameter standing in for $\epsilon$.

Then, with $\epsilon > 0$, by splitting up the area under the curve into two regions, one to the left and one to the right of the vertical line $x=x_1$, and shifting variables $x=x_1+t$, we get
\begin{eqnarray*}
    F(x_1 + \epsilon) &=& 
    \int_{x_0}^{x_1} f(x) \, dx + \int_{x_1}^{x_1 + \epsilon} f(x) \, dx \mbox{, splitting the region in two} \\
    &=& F(x_1) + \int_0^\epsilon f(x_1+t)\, dt 
    \mbox{, shifting the second integral} \\
    &=& F(x_1) + \int_0^\epsilon f(x_1)+E(t)\, dt \mbox{, using the error function}\\
    &=& F(x_1) + \epsilon \cdot f(x_1) + \int_0^\epsilon E(t)\, dt
    \mbox{, integrating the constant term.} 
\end{eqnarray*}
So we have obtained a linear approximation for $F(x_1 + \epsilon)$ with scaling factor $f(x_1)$, which shows the derivative of $F$ at $x_1$ is exactly $f(x_1)$. That is, $F'(x_1) = f(x_1)$, which is the statement of the FTC that we wanted to show. 

Of course, we must still show the last term $\int_0^\epsilon E(t)\, dt$ is of the form $|\epsilon| E_1(\epsilon)$ for some error function $E_1(\epsilon).$ But that is pretty easy. Define two functions $m(\epsilon) = \min_{ \{ 0\leq t \leq \epsilon \} } E(t)$ and  $M(\epsilon) = \max_{ \{ 0\leq t \leq \epsilon \}} E(t)$. Then by integrating the inequalities $m(\epsilon) \leq E(t) \leq M(\epsilon)$ we get
\[ m(\epsilon) \epsilon \leq \int_0^\epsilon E(t)\, dt \leq M(\epsilon)\epsilon.\]
Thus we define $E_1(\epsilon) = \frac{1}{\epsilon} \int_0^\epsilon E(t)\, dt$ and we have
\[ m(\epsilon) \leq E_1(\epsilon) \leq M(\epsilon).\]
The fact that $E(\epsilon)$ is an error function implies its min and max functions $m(\epsilon), M(\epsilon)$ are error functions and thus so is $E_1(\epsilon). $
This completes the demonstration.

\section{Concluding remarks}

Our goal here has been to address the serious concern that students in first year calculus find limits both a source of difficulty and a distraction from the main ideas and themes of calculus. We have proposed a method focused on a precise definition of approximation and demonstrated its utility with a series of examples from the major topics in calculus that are rigourous and yet also natural and intuitive. We feel this approach will alleviate some of the difficulties the students are having with limits and allow them to master the fundamentals of calculus without watering down the content in any way. 

Along the way, we have discovered consistent and elegant approaches to defining the trigonometric, exponential and hyperbolic functions from simple geometry via an intuitive use of areas under a curve. This approach also produces a unified technique for describing the properties of these -- at first glance -- dissimilar but basic functions of calculus. 

Our belief is that this foundational method of approximation can be expanded to a complete course in calculus, including differential and integral calculus. It would be interesting to explore its use in courses in differential equations and modern analysis. 

We would like to thank the many colleagues who encouraged us to put these ideas down in hard copy, and especially Professor Emeritus John Fournier of the University of British Columbia who directed us to earlier works in ``limitless calculus.''

\bibliographystyle{plain} 
\bibliography{references}

\end{document}